\documentclass{amsart}

\usepackage[english]{babel}

\usepackage[letterpaper,top=2cm,bottom=2cm,left=3cm,right=3cm,marginparwidth=1.75cm]{geometry}

\usepackage{amsmath}
\usepackage{graphicx}
\usepackage[colorlinks=true, allcolors=blue]{hyperref}
\usepackage[english]{babel}
\usepackage{amssymb}
\usepackage{xcolor}
\usepackage{bbm}
\usepackage{bbold}
\usepackage{xcolor}  
\usepackage{mathrsfs}

\newcommand{\R}{\mathbb{R}^3}
\newcommand{\RS}{\mathbb{R}^3\times\mathbb{S}^2}
\newcommand{\B}{\mathcal{B}}
\newcommand{\E}{\mathcal{E}}
\newcommand{\M}{\mu}
\newcommand{\dd}{\mathrm{d}}
\newcommand{\p}{\partial}
\newcommand{\hx}{\hat{x}}
\newcommand{\hv}{\hat{v}}
\newcommand{\hw}{\hat{w}}
\newcommand{\n}{\sigma}
\newcommand{\hf}{\hat{f}}
\newcommand{\hF}{\hat{F}}
\newcommand{\id}{\mathrm{Id}}
\newcommand{\Div}{\mathrm{div}}
\newcommand{\q}{q\left(\frac{v-u}{\sqrt{\theta/m_g}}\right)}
\newcommand{\fed}{f^{\eta,\delta}}
\newcommand{\Fed}{F^{\eta,\delta}}
\newcommand{\num}{n}

\newtheorem{theorem}{Theorem}

\newtheorem{proposition}{Proposition}

\title{A Formal Derivation of the Thick Spray Model from the Enskog-Boltzmann System}
\author{Zhe Chen}
\address{Sorbonne Universit\'e, LJLL, Bo\^\i te courrier 187, 75252 Paris Cedex 05 France}
\email{zhe.chen@sorbonne-universite.fr}

\begin{document}


\maketitle
\begin{abstract}
We propose a specific scaling that formally derives the Euler-Vlasov model for thick sprays which is widely adopted in engineering from the Boltzmann-Enskog model. Beyond validating the kinetic-fluid equations underlying this model, we also identify higher-order corrections of the viscous force (of the order of the volume of the dispersed phase), which is neglected due to the incorporation of the volume fraction (void) function. Our approach builds on the work [Laurent Desvillettes, Fran\c{c}ois Golse, and Valeria Ricci. A formal passage from a system of Boltzmann equations for mixtures towards a Vlasov-Euler system of compressible fluids. Acta Mathematicae Applicatae Sinica, English Series, 35(1):158–173, Jan 2019.], who formally connected a system of coupled Boltzmann equations for  binary mixtures to a Vlasov-Euler system for thin sprays.
\end{abstract}

\section*{Introduction}

An aerosol, or spray, is a binary mixture composed of a dispersed phase--- typically liquid droplets, or solid particles---suspended in a gaseous medium known as the propellant. There are various types of equations to describe such a system; see, for instance \cite{desvillettes2010}. A significant class of models describing the dynamics of aerosol or spray flows couples
\begin{enumerate}
    \item[(a)] a kinetic equation governing the dispersed phase, with
    \item[(b)] a fluid equation describing the evolution of the propellant.
\end{enumerate}
These fluid-kinetic equations, also known as Eulerian-Lagrangian or gas-particles model, are coupled through the drag force. 

According to the proportions of the dispersed phase, sprays regimes can be classified into three categories (see, for instance, Chapter 1 in \cite{ORourke1981}): the very thin sprays, where the volume fraction of the dispersed phase's proportion is much less than $10^{-3}$; the thin sprays, where the proportion is much less than $10^{-1}$; and thick sprays, where the dispersed phase occupies up to approximately $10^{-1}$ of the total volume.
From a microscopic point of view, let $N_p$ denotes the total number of particles and $r_p$ the radius of each particle, these regimes correspond to different scaling behaviors: very dilute scale $r_pN_p=O(1)$; the Boltzmann-Grad scale $r_p^2N_p=O(1)$; dense gas scale $r_p^3N_p=O(1)$. Any dense scaling, such as, $N_p r_p^4=O(1)$, would imply that the total volume of the dispersed phase, $N_p r_p^3$, diverges as $r_p \to 0$, which is physically inconsistent.

Intuitively speaking, as the concentration of the dispersed phase (dusts or droplets) increases, the interaction between two species becomes stronger. 

In the case of a \emph{very thin spray}, we only take into account interactions of order $O(r_p)$. This, for instance, gives rise to the Brinkman force; see \cite{Allaire1991,DGR2008}. However, in the Euler regime, the friction force appears only at order $O(r_p^2)$. As a result, the gas phase dynamics can be computed independently of the influence of the droplets.

When the spray becomes \emph{thin}, the proportion of the dispersed phase grows, and its feedback on the gas phase become non-negligible---although its volume remains small. In this intermediate regime, the dispersed phase influences the exchange rates of momentum and energy with gas. The viscous force emerges at this stage and is intuitively proportional to the surface area of the spherical particles, i.e. $O(r_p^2)$. At this level of approximation, the volume of the dispersed phase, which scales as $O(r_p^3)$, is still neglected. \cite{DBGR2017,DBGR2018} provided the formal derivation from coupled Boltzmann equations to Vlasov-Navier-Stokes, and Vlasov-Stokes respectively. \cite{DGR2019} provided the formal derivation from the same system to Vlasov-Euler equations.

In the \emph{thick sprays} regime, the volume of particles become significant. Consequently, the drag force between two species becomes more complex---it includes not only the viscous effects but also the modification of the pressure due to the volume occupation. This leads to the the notion of ``boundary'' between dispersed phase and the propellant. On such boundaries, pressure imbalance may arise in non-equilibrium state.
To address this, \cite{THA1976} introduced the concept of the volume fraction (or void) function, which behaves like $1-O(r_p^3)$. This modifies the drag force to include not only viscous effects but also corrections associated with volume exclusion.

A key contribution of this work is the derivation of higher-order perturbations (of order $O(r_p^3)$) in the viscous force from a mesoscopic perspective—corrections that were missed in previous models such as those proposed in \cite{ORourke1981,DUKOWICZ1980}. Moreover, we provide a justification for the plausibility of these extended models from the kinetic–kinetic framework.

   The outline of this paper is as follows: we introduce the equation of Boltzmann and its variant---the Enskog equation---are presented in Section \ref{sec:introduction to boltzmann equation} and \ref{sec:introduction to enskog equation}, respectively. The Section \ref{sec:kinetic-fluid equation} introduce the target equations which we aim to derive. In Section \ref{sec:building on the dimensionless kinetic equations}, we formulate the initial coupled kinetic system. In section \ref{sec:Main result}, we state the main result of this paper. Section \ref{sec:Particle equation} demonstrates how to derive the Vlasov equation from the Enskog-Boltzmann equation. In section \ref{sec:Gas phase equations}, we establish the conservation laws for the gas phase using the moment method. Finally, we conclude the result with some remarks in Section \ref{sec:conclusion}.



\section{The Boltzmann equation and Enskog-Boltzmann equation}

Let $f\equiv f(t,x,w)\ge 0$ and $F\equiv F(t,x,v)\ge 0$ denote the number distribution functions of gas molecules and particles, respectively, depending on time $t\in \mathbb{R}_+$, position $x\in \R$, and velocities $w,v\in\R$. 

Calligraphic letters such as $\B$ and $\E$ are used to denote the collision integrals. In the notation used throughout this paper, square brackets \([\,]\) indicate functional dependence, while parentheses \((\,)\) indicate dependence on variables.

\subsection{The Boltzmann equation}\label{sec:introduction to boltzmann equation}
\begin{figure}
    \centering
    \includegraphics[width=0.6\textwidth]{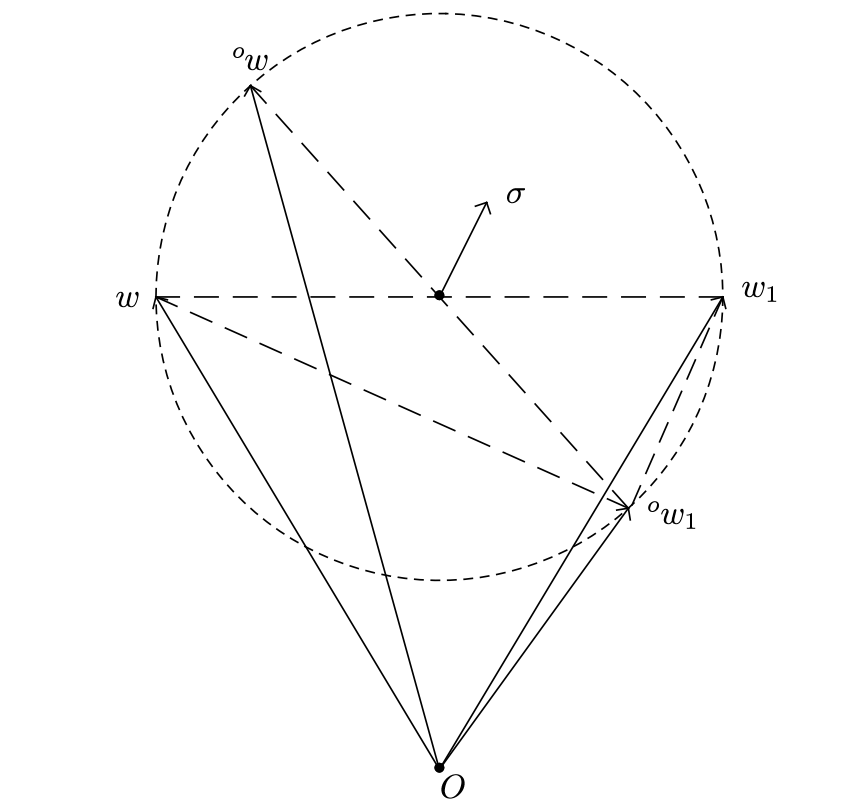}
    \caption{Illustration of a geometric configuration relevant to Boltzmann collisions}
    \label{fig:boltzmann_collision}
\end{figure}

The Boltzmann equation is the classical collisional kinetic model describing the dynamics of a dilute gas, and it takes the following form:
\begin{equation}\label{the boltzmann equation}
\p_t f(t,x,w)+w\cdot\nabla_xf(t,x,w)=\mathcal{B}[f,f](t,x,w),
\end{equation}
where the Boltzmann collision integral is 
\begin{equation}\label{original boltzmann collision integral}
\mathcal{B}[f,f](w)=\iint_{\RS}\bigg(f({}^ow)f({}^ow_1)-f(w)f(w_1)\bigg)b(w-w_1, \n)\dd w_1\dd \n.
\end{equation}
Here ${}^ow$ and ${}^ow_1$ denote the velocities just before the collision, while $w$ and $w_1$ are the velocities immediately after the collision. The positive term $f({}^ow)f({}^ow_1)$ in the collision integral above accounts for the ``gain'' of gas molecules with velocity $w$, while the negative term $-f(w)f(w_1)$ represents the ``lose'' of such molecules with velocity $w$. For elastic collision between gas molecules, it's natural to require conservation of momentum and energy:
\[
{}^ow+{}^ow_1=w+w_1;\qquad |{}^ow|^2+|{}^ow_1|^2=|w|^2+|w_1|^2.
\]
The system involves $6$ unknowns ($3$ component for each pre-collision velocity), and $4$ constraints, leaving $2$-degree of freedom. To parametrize ${}^ow$ and ${}^ow_1$, we introduce additional variable $\n\in\mathbb{S}^2$, and it's easy to verify the following representation satisfy the conservation of momentum and energy:
\begin{align*}
    &{}^ow\equiv{}^ow(w,w_1,\n)=w-(w-w_1)\cdot \n \n,\\
    &{}^ow_1\equiv{}^ow_1(w,w_1,\n)=w_1+(w-w_1)\cdot \n \n,
\end{align*}
where $\n$ varies over the unit sphere $\mathbb{S}^2 \subset \mathbb{R}^3$.
Note that here the pre-collision velocities are parametrized by the unit vector $\n$. However, this is not the standard $\sigma$-representation commonly used in the Boltzmann collision integral.
The representation of the pre-collision velocities ${}^ow$ and ${}^ow_1$ are not unique. The relationship between these vectors can be visualized in Figure \ref{fig:boltzmann_collision}. 

The collision kernel $b(z,\omega)$ is of the form
\[
b(z,\omega)=|z|\Sigma(|z|,\cos(\widehat{z,\omega}))\qquad \mathrm{with} ~\Sigma>0.
\]
One can (formally) derive the Boltzmann equation from Hamilton's equation for an $N$-particle system in the case of hard sphere---rigorously established for the short time by Lanford \cite{Lanford1975}---under the Boltzmann-Grad scaling: $N(2r)^2\sim 1$. In this case, collision kernel takes of the form $b(z,\omega)=(2r)^2|\cos(\widehat{z,\omega})|$. Moreover, the excluded volume $N(2r)^3=O(r)$ becomes negligible in the $N$-particle limit (namely, $N\to \infty$ and $r\to 0$). This indicates that the gas described by the Boltzmann equation must be sufficiently ``dilute''.

It's well-known that for all $f\in L^1(\R,(1+|w|^2)^3\dd w)$, the collision integral conserves mass, momentum, energy locally (see, for instance, section 3.1 and section 3.3 in \cite{cercignani1994}), meaning
\begin{equation}
    \int_{\R}\B[f,f](w)\begin{pmatrix}1\\w\\|w|^2\end{pmatrix}\dd w=0.
\end{equation}
This property implies that the Boltzmann collision integral does not contribute to the macroscopic equation \emph{directly}, but only in an indirect way. In other words, we multiply equation \eqref{the boltzmann equation} by $1$, $w$, or $|w|^2/2$ and integrate it with respect to $w$. If $f,\p_t f$ and $\nabla_xf\in L^1(\R,(1+|w|^2)^3\dd w) $  we can write

\begin{equation*}\left\{
    \begin{aligned}
        &\p_t\int_{\R} f(t,x,w)\dd w+\Div_x\int_{\R} wf(t,x,w)\dd w=0,\\
        &\p_t\int_{\R} wf(t,x,w)\dd w+\Div_x\int_{\R} w^{\otimes 2}f(t,x,w)\dd w=0,\\
        &\p_t\int_{\R} \frac{|w|^2}{2}f(t,x,w)\dd w+\Div_x\int_{\R} w\frac{|w|^2}{2}f(t,x,w)\dd w=0,
    \end{aligned}\right.
\end{equation*}
where the divergence of a $2$-tensor $A\equiv(A_{ij})$ with respect to variable $x$ is defined componentwise by $(\Div_x A)_{i}=(\p_{x_j} A_{ij})$.

Under the following definitions, we can establish the connection between kinetic description of gas and the framework of continuum mechanics. Let's introduce:
\begin{equation*}
    \begin{aligned}
        &R(t,x):=\int_{\R} f(t,x,w)\dd w,\qquad &\mathrm{Local ~number ~density}\\
        &U(t,x):=\frac{H(R(t,x))}{R(t,x)}\int_{\R} wf(t,x,w)\dd w,\qquad &\mathrm{Bulk ~velocity}\\
        &P(t,x):=\int_{\R} (w-U(t,x))^{\otimes 2}f(t,x,w)\dd w,\qquad &\mathrm{Stress ~tensor}\\
        &E_{\mathrm{int}}(t,x):=\frac{H(R(t,x))}{2R(t,x)}\int_{\R}|w-U(t,x)|^2f(t,x,w)\dd w,\qquad &\mathrm{Internal ~energy}\\
        &\Lambda(t,x):=\frac{1}{2}\int_{\R}(w-U(t,x))|w-U(t,x)|^2f(t,x,w)\dd w,\qquad &\mathrm{Energy ~flux}
    \end{aligned}
\end{equation*}
where $H:\mathbb{R}\to\{0,1\}$ denotes the Heaviside step function is defined by
\[
H(r)=\left\{\begin{aligned}
     &1 &\mbox{for}~ &r>0  \\
     &0 &\mbox{for}~ &r\le 0.
\end{aligned}\right.
\] 
These definitions leads to the following fluid equations:
\[
\left\{
    \begin{aligned}
        &\p_t R + \Div_x(RU) = 0,\\
        &\p_t(RU_i)+\sum_{j=1}^3\p_{x_j}(RU_iU_j+P_{ij})=0,\qquad i=1,2,3\\
        &\p_t\left(R\left(\tfrac{1}{2}|U|^2 + E_{\mathrm{int}}\right)\right) + \Div_x\left(RU\left(\tfrac{1}{2}|U|^2 + E_{\mathrm{int}}\right)\right)+\sum_{i=1}^3 \p_{x_i}\left[\sum_{j=1}^3U_jP_{ij}+ \Lambda_i\right] = 0,
    \end{aligned}
\right.
\]
These equations represent, respectively, the conservation of mass, momentum, and energy in the fluid regime derived from kinetic theory. Obviously, the system of macroscopic equation is not closed, as it involves $16$ unknowns subject to only $5$ constrains. To close the system, we consider the high-collision regime, which ensures that $f$ take the form of a local Maxwellian distribution (see, for instance, Section 3.2 in \cite{cercignani1994}). Under this assumption, it's straightforward to verify the following relations:
\[
P(t,x)=\frac{2}{3}R(t,x)E_{\mathrm{int}}(t,x)\id;\qquad \Lambda(t,x)=0.
\]
These constitute the classical equation of state for dilute (or ideal) gas.

This illustrates how the Boltzmann collision integral contributes to the fluid equations, in contrast to the free transport equation $\p_tf+w\cdot\nabla_xf=0$, where no such closure naturally arises---yet the stress tensor still appears in the macroscopic description.

For the binary mixture, we can consider the following coupled Boltzmann equations:
\begin{equation}\label{Original coupled Boltzmann equations}
\left\{
\begin{array}{ll}
    (\partial_t + w \cdot \nabla_x) f (t, x, w) &= \mathcal{R}[f,F](t,x,w)+\mathcal{C}[f,f](t,x,w),\\
    (\partial_t + v \cdot \nabla_x) F (t, x, v) &= \mathcal{D}[F,f](t,x,v)+\mathcal{Q}[F,F](t,x,v).
\end{array}
\right.
\end{equation}
   Here, the collision integrals $\mathcal{R},\mathcal{D}$ describe the interactions between two different species, while collision integrals $\mathcal{Q},\mathcal{C}$ represent the self-collision within each species. 
   Previous works are largely based on specific choices and modifications on the Boltzmann collision integral. For example, 
   \cite{DBGR2017} and \cite{DBGR2018} use the Boltzmann collision integrals with different scalings for velocities and mass ratio.
   \cite{CD2009} considers the size of particles in the modeling.
   \cite{CD2024} and \cite{CMS2025} modified the Boltzmann collision integral. In this work, we adopt the Enskog collision integral, which will be introduced in the next subsection, as a refinement of the Boltzmann operator to more accurately capture the finite-volume effects in moderately to densely packed regimes.

\subsection{The Enskog equation}\label{sec:introduction to enskog equation}

When the radius of molecules cannot be ignored, precisely, the Boltzmann-Grad assumption fails, we should seek for another model to describe such a system (dense gas). Later on, Enskog \cite{enskog1922} proposed the so-called (standard) Enskog equation which shares the same structure with the Boltzmann equation, while the collision integral is different:
\begin{equation}\label{Enskog equation}
    (\p_t+v\cdot\nabla_x)f(t,x,v)=\E[f,f](t,x,v),
\end{equation}
where collision integral $\E$ is given by
\begin{multline}\label{y representation of enskog collision integral}
    \E[f,g](t,x,v):=\int_{\RS}\bigg[\chi[f,g](t,x,y) f(t,x,{}^ov)g(t,y,{}^ov_1)-\chi[f,g](t,x,y) f(t,x,v)g(t,y,v_1)\bigg]\\
    (v-v_1)\cdot \frac{y-x}{|y-x|}H\left((v-v_1)\cdot \frac{y-x}{|y-x|}\right)\delta(|x-y|-2r_p)\dd y \dd v_1.
\end{multline}
Here, $r_p > 0$ denotes the diameter of the particle, and $\chi \equiv \chi[f, g](t, x, y)$ is the correlation function, which captures the effects of volume exclusion—such as shielding influence. Further details can be found in Chapter 16 of \cite{chapman1990}.
 In other words, the presence of finite-size particles reduces the statistical independence of two-particle configurations compared to the dilute gas case.
If we assume the propagation of chaos holds---that is, the two-particle distribution function $f^{(2)}(t,x_1,x_2,w_1,w_2)$ can be approximated by the product of two single-particle distribution functions $f(t,x_1,w_1)f(t,x_2,w_2)$---then we may set $\chi\equiv 1$ in the Enskog collision integral. 
Under this assumption, we refer to equation \eqref{Enskog equation} with such a simplified collision integral as the \emph{Enskog-Boltzmann equation}. For the remainder of this work, we consider only this case, i.e. we always assume $\chi\equiv 1$. As we will see, this assumption is sufficient to capture  the necessary modification of the pressure in the macroscopic fluid equations.

Unlike in the Boltzmann case, the pre-collision velocities in the Enskog setting are determined uniquely by conservation momentum, conservation of energy, and the positions of the colliding particles
\begin{equation*}
\left\{
  \begin{aligned}
    &{}^ov\equiv {}^ov(x,v_1,x-y)=v-(v-v_1)\cdot\frac{y-x}{|y-x|}\frac{y-x}{|y-x|},\\
    &{}^ov_1\equiv {}^ov_1(x,v_1,x-y)=v_1+(v-v_1)\cdot\frac{y-x}{|y-x|}\frac{y-x}{|y-x|}.
\end{aligned} \right. 
\end{equation*}
\begin{figure}
    \centering
    \includegraphics[width=0.65\textwidth]{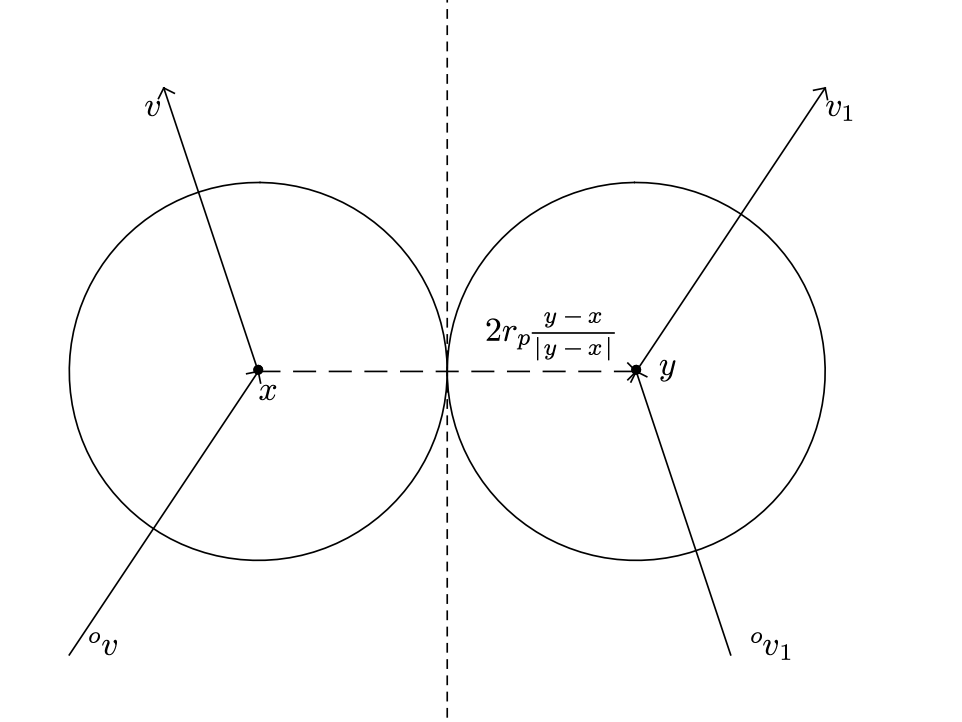}
    \caption{Diagram illustrating particle positions and velocities during a collision.}
    \label{fig:particle_collision}
\end{figure}
This determinism arises because, unlike in the Boltzmann model where particles are treated as point masses, the Enskog model accounts for the finite size of particles. Thus, it distinguishes between the collision point and the center of mass of each particles. See the illustrative diagram \ref{fig:particle_collision}. More details about Enskog equation can be found in Chapter 16 of \cite{chapman1990}. 
In the following, we adopt the $\n$-representation of the Enskog-Boltzmann collision integral, which corresponds to the change of variables $y = x \pm 2r_p \n$ in the gain (positive) and loss (negative) terms, respectively. Under this transformation, the collision integral \eqref{y representation of enskog collision integral} can be rewritten as:
\begin{multline*}
    \E[f,g](t,x,v):=\int_{\RS}\bigg[ f(t,x,{}^ov)g(t,x-2r_p \n,{}^ov_1)- f(t,x,v)g(t,x+2r_p \n,v_1)\bigg]\\
    (v-v_1)\cdot \n H\left((v-v_1)\cdot \n\right)(2r_p)^2\dd \n \dd v_1.
\end{multline*}

Parallel to the result for the Boltzmann collision, integral, one may ask whether the Enskog collision integral locally conserves mass, momentum, and energy? The answer is \emph{no}. However, it does conserve these quantities globally.  As shown in \cite{CCG2024}, the Enskog-Boltzmann collision integral admits the following conservative form (in the sense of distributions):
\[
\E[f,f](w)\begin{pmatrix}1\\w_1\\w_2\\w_3\\|w|^2\end{pmatrix}=\begin{pmatrix}\Div_v\mathbb{J}_0\\\Div_x\mathbb{I}_1+\Div_v\mathbb{J}_1\\\Div_x\mathbb{I}_2+\Div_v\mathbb{J}_2\\\Div_x\mathbb{I}_3+\Div_v\mathbb{J}_3\\\Div_x\mathbb{I}_4+\Div_v\mathbb{J}_4\end{pmatrix},
\]
for some specific vectors $\mathbb{I}_1,\dots,\mathbb{I}_4,\mathbb{J}_0,\dots,\mathbb{J}_4$.
This indicates that, unlike the Boltzmann case, the Enskog collision integral contribute \emph{directly} to the macroscopic equations due to the presence of the term $\mathbb{I}$.
More specifically, using the same notations as before, we obtain the following macroscopic equations by multiplying equation \eqref{Enskog equation} by $1,w$, or $|w|^2/2$, and integrate it with respect to $w$ 
\[
\left\{
    \begin{aligned}
        &\p_t R + \Div_x(RU) = 0,\\
        &\p_t(RU_i)+\sum_{j=1}^3\p_{x_j}(RU_iU_j+P_{ij})-\Div_x\mathbb{I}_i=0,\qquad i=1,2,3\\
        &\p_t\left(R\left(\tfrac{1}{2}|U|^2 + E_{\mathrm{int}}\right)\right) + \Div_x\left(RU\left(\tfrac{1}{2}|U|^2 + E_{\mathrm{int}}\right)\right)+\sum_{i=1}^3 \p_{x_i}\left[\sum_{j=1}^3U_jP_{ij}+ \Lambda_i\right]-\Div_x\mathbb{I}_4 = 0,
    \end{aligned}
\right.
\]
Another approach to derive the corresponding hydrodynamics equations for Enskog gas---Enskog-Euler system---is through the Hilbert expansion, as in \cite{Lachowicz1998}. However, unlike the system derived above, this approach yields an asymptotic result.
As a byproduct, the influence of particle volume on the equation for the Enskog gases becomes apparent. The pressure is modified as follows
\begin{equation}\label{EOS for Enskog gas}
p_E=\frac{2}{3}RE_{\mathrm{int}}\left(1+\frac{4}{3}\pi r_p R \right),
\end{equation}
highlighting a correction to the ideal gas law due to finite-size effects. At this point, we have strong motivation to introduce the Enskog collision integral to capture the effect of volume fraction on stress tensor in the macroscopic equation for thick sprays.


\section{The kinetic-fluid description: Vlasov-Euler system}\label{sec:kinetic-fluid equation}
Beyond the kinetic description of collisions, an alternative framework involves coupling a fluid equation---such as the Navier-Stokes, Stokes, or Euler equations---with a kinetic equation, typically of Vlasov type. When the volume fraction of the dispersed phase is sufficiently small, the system is described by the so-called thin sprays model:

\begin{equation}\label{Vlasov-Euler for thin sprays}
   \left\{
   \begin{aligned}
       &\p_t \rho+\Div_x(\rho u)=0,\\
       &\p_t (\rho u)+\Div_x(\rho u^{\otimes 2})+\nabla_x p=m_g\int_{\R}F(v)D(v-u)\dd v,\\
       &\p_t\left( \rho E\right)+\Div_x\left( \rho uE\right)+\Div_x( u p)=m_g\int_{\R}v\cdot D(v-u) F(v)\dd v,\\
       &\p_t F+v\cdot\nabla_x F+\Div_v[F(v) D(v-u)]=0,\\
       &p=\num\theta,\,E=\frac{|u|^2}{2}+\frac{3\theta}{2m_g},\, n=\frac{\rho}{m_g}
   \end{aligned}\right.
   \end{equation}
where $\rho,u,\theta,p,E$ represent the local density, velocity, temperature, pressure, and internal energy of gas respectively, $F=F(t,x,v)$ denotes the distribution of particles bedding in the gas, $D\equiv D(v-u)$ is a given function describing the friction between the two species, typically defined by the linear case $D(\xi)=\xi$, or the quadratic case $D(\xi)=\xi|\xi|$. However, it may depending on gas density and internal energy. A formal derivation of this model has been presented in \cite{DGR2019}.

To account for the volume contribution in the pressure, we introduce the volume fraction function by 
\[
\check{\alpha}(t,x):=1-\frac{4\pi}{3}r_p^3\int_{\R}F(t,x,v)\dd v,
\]
where $r_p$ represent the radius of spherical suspension, and $F(t,x,v)$ is the density function of suspensions. Thus, $\int_{\R}F(t,x,v)\dd v$ represents the total number of suspensions. Here, we give an intuitive explanation of how the volume fraction contributes when it cannot be neglected. From a macroscopic perspective, one can understand why volume occupation contributes to the pressure; see \cite{THA1976} for reference. The buoyancy force on a spherical droplet can be expressed as
    $(F_{\mathrm{buoy}})_i=-\int_{\p B_i} p(t,x)n(x)\dd S(x)$, 
    where $\dd S$ is surface element, $\n$ denotes the outer normal vector to $i$-th ball $B_i=B(X_i,r_p)$. By the Gauss theorem, we obtain
    \((F_{\mathrm{buoy}})_i=-\int_{B_i}\nabla_{x}p(t,x)\dd V(x).\)
    If the density of the droplet does not vary significantly on its surface, or if the droplet is sufficiently small, the buoyancy can be approximated by
    \(
    (F_{\mathrm{buoy}})_i\sim -\frac{4\pi}{3}r_p^3\nabla_xp.
    \)
    Repeating this argument for $N_p$ particles, we obtain
    \[
    F_{\mathrm{buoy}}\sim -\frac{4\pi}{3}r_p^3N_p\nabla_xp.
    \]
    From the mesoscopic point of view, the number of particles is replaced by the particle distribution function $F$, leading to the approximation  
    \[
    F_{\mathrm{buoy}} \sim -\frac{4\pi}{3}r_p^3 \int_{\R} F(t,x,v) \, \mathrm{d}v \, \nabla_x p.
    \]
   This explains, in an extremely formal way, why the volume fraction function appears outside the pressure gradient in the following system. We will revisit and discuss this reasoning from a mesoscopic point of view in Section \ref{sec:conclusion}. Therefore, we can formally write the following Vlasov-Euler equations for thick sprays which is widely accepted:
   \begin{equation}\label{Vlasov-Euler for thick sprays}
   \left\{
   \begin{aligned}
       &\p_t(\check{\alpha} \rho)+\Div_x(\check{\alpha}\rho u)=0,\\
       &\p_t(\check{\alpha} \rho u)+\Div_x(\check{\alpha}\rho u^{\otimes 2})+\nabla_x p=m_g\int_{\R}F(v)\Gamma(v-u)\dd v,\\
       &\p_t\left(\check{\alpha} \rho E\right)+\Div_x\left(\check{\alpha} \rho uE\right)+\Div_x(\check{\alpha} u p)+p\p_t\check{\alpha}=m_g\int_{\R}D(v-u)\cdot vF(v)\dd v,\\
       &\p_t F+v\cdot\nabla_x F+\Div_v(F \Gamma)=0,\\
       &\Gamma(v-u)=\frac{4\pi}{3}r_p^3\nabla_xp+D(v-u),\\
       &p=n\theta,\,E=\frac{|u|^2}{2}+\frac{3\theta}{2m_g},\, \rho=m_g n,\\
       &\check{\alpha}=1-\frac{4\pi}{3}r_p^3\int_{\R}F(v)\dd v,
   \end{aligned}\right.
   \end{equation}
   where $\Gamma$ represents the drag force between gas and particle. Here the drag force includes not only the friction term $D$, as in the case of thin sprays, but also an additional contribution accounting for pressure modification.
   If we set $\check{\alpha}=1$ in equations \eqref{Vlasov-Euler for thick sprays}, we recover the Vlasov-Euler system for thin sprays \eqref{Vlasov-Euler for thin sprays}. 
  
   A substantial body of work has focused on modeling and analysis of gas-particle systems, particularly thick sprays.
   The foundational model of many particles suspended in a gas was considered in \cite{THA1976}. A more systematic and comprehensive model for thick sprays---including collisions, cooling process, energy exchange, and more--- was developed in \cite{ORourke1981}. This model has served as the basis for many subsequent studies and computational developments, especially in the content of the Kiva code; see, for instance \cite{AOB1989}.
    \cite{BDM2003} conducted a numerical study of this model, and demonstrated the global conservation of total energy. \cite{BDGN2012} provided a a one-dimensional simulation of a thick sprays system within a pipe.
  
   From a theoretical perspective, \cite{Fournet2024} established the well-posedness for the modified thick sprays model by averaging the pressure term and volume fraction. \cite{EH2023} shows the well-posedness of this system under initial data satisfying a Penrose-type stability, which is shown to be necessary for the long-time stability. \cite{BDF2024,fournet2024b} studied the linearized thick sprays model, they showed that the system have similar behavior as Landau damping, in both numerical and theoretical ways. More recently, \cite{BDFG2025} shows the thick sprays model can exhibit linear ill-posedness, even locally in time.

 Beyond this system, \cite{DM2010} formally derive the Eulerian-Eulerian description for the thick sprays.

 Compared to general gas-particle models without finite-volume effect, the setting we consider involves the following simplifications:  the radius of particles is fixed; there isn't energy exchanged, and collision between particles are elastic; all the collisions are elastic.The case with variable radii is addressed in \cite{CD2009}. A model incorporating inelastic collisions can be found in \cite{CDS2012} with numerical scheme. Models that include energy exchange between different species have been studied in more recent works, such as \cite{CD2024,CMS2025}.
   
\section{Building on the dimensionless kinetic equations}\label{sec:building on the dimensionless kinetic equations}
 As motivated in the last subsection, to model the thick sprays, we introduce the Enskog-Boltzmann collision integral into the system \eqref{original boltzmann collision integral} with the following (elastic) collision laws:
   \begin{align}
        \label{collision conserves momentum}
    & m_pv' + m_g w' =  m_pv + m_g w\\
     \label{collision conserves energy}
    & m_p|v'  |^2+ m_g |w'|^2 = m_p|v|^2 + m_g |w|^2,
    \end{align}
Where $m_p$ and $m_g$ denote the mass of particle and gas molecule, respectively.
These velocity distribution functions should satisfy the coupled Enskog-Boltzmann equations:
\begin{equation}\label{Original coupled Enskog-Boltzmann equations}
\left\{
\begin{array}{ll}
    (\partial_t + w \cdot \nabla_x) f (t, x, w) &= \E_1[f,F](t,x,w)+\E_0[f,f](t,x,w)\\
    (\partial_t + v \cdot \nabla_x) F (t, x, v) &= \E_2[F,f](t,x,v)+\E_3[F,F](t,x,v)
\end{array}
\right.
\end{equation}
where the Enskog-Boltzmann collision integrals are given by
\begin{equation*}\label{Boltzmann collision integral for gas}
    \E_0[f,f](x,w)=\int_{\RS} [f (x, {}{}^ow) f (x-2r_g\n, {}{}^ow_1) - f (x , w) f (x+2r_g\n, w_1)]   b_{g}( w - w_1, \n) \dd  \n \dd  w_1,
\end{equation*}
\begin{multline*}
    \E_1[f,F](x,w)=(r_g+r_p)^2\int_{\RS} [F (x + (r_g+r_p) \n, {}'v) f (x, {}'w) - F (x - (r_g+r_p) \n, v) f (x, w)]\\
    (v-w)\cdot \n  H(( v - w) \cdot \n) \dd  \n \dd  v
\end{multline*}
\begin{multline*}\label{Enskog-Boltzmann collision integral for particle}
\E_2[F,f](x,v)=(r_g+r_p)^2\int_{\RS} [f (x - (r_g+r_p) \n, {}'w) F (x, {}'v) - f (x + (r_g+r_p) \n, w) F (x, v)]\\
(v-w)\cdot \n  H(( v - w) \cdot \n) \dd  \n \dd  w
\end{multline*}
\begin{equation*}\label{Boltzmann collision integral for particle}
    \E_3[F,F](x,v)=\int_{\RS}  [F(x, {}{}^ov) F (x-2r_p\n, {}{}^ov_1) - F (x , v) F(x+2r_p\n, v_1)] b_{p}( v - v_1, \n) \dd  \n \dd  v_1,
\end{equation*}
where the collision kernel $b_{\beta}$ is of the form $b_{\beta}(\xi,\sigma)=|\xi|\Sigma_\beta(|\xi|,\sigma).$
For the case of hard sphere, $\Sigma_{\beta}(\xi,\sigma)$ is given by $\frac{2r_{\beta}^2}{m_\beta}|\cos (\widehat{\xi,\sigma})|$. Here $r_{\beta}$ denotes the radius of molecule. Specifically, $\beta=g$ and $p$ correspond to the collision kernel $b_g$ and $b_p$ respectively, and pre-collosion velocities between same species are given by
\[ \left\{ \begin{array}{l}
     {}{}^ov={}{}^ov (v, v_1, \n) = v - \left( v -  v_1\right) \cdot \n \n\\
     {}{}^ov_1={}{}^ov_1(v, v_1, \n) = v_1 -  (v_1-  v) \cdot \n \n,   \end{array} \right. \]
 and the pre-collision velocities between different species are given by
 \[
 \left\{ \begin{array}{l}
     {}'v={} 'v(v, w,\n, \eta) = v - \frac{2\eta}{1+\eta}\left( v -  w\right) \cdot \n \n\\
     {}'w={}'w(v, w, \n,\eta) = w -  \frac{2}{1+\eta}(w-  v) \cdot \n \n ,  \end{array} \right.
 \]
 where $\eta=\frac{m_g}{m_p}$ denotes the mass ratio, supposing to be sufficient small.
 
 \begin{figure}[htbp]
    \centering
    \includegraphics[width=0.7\textwidth]{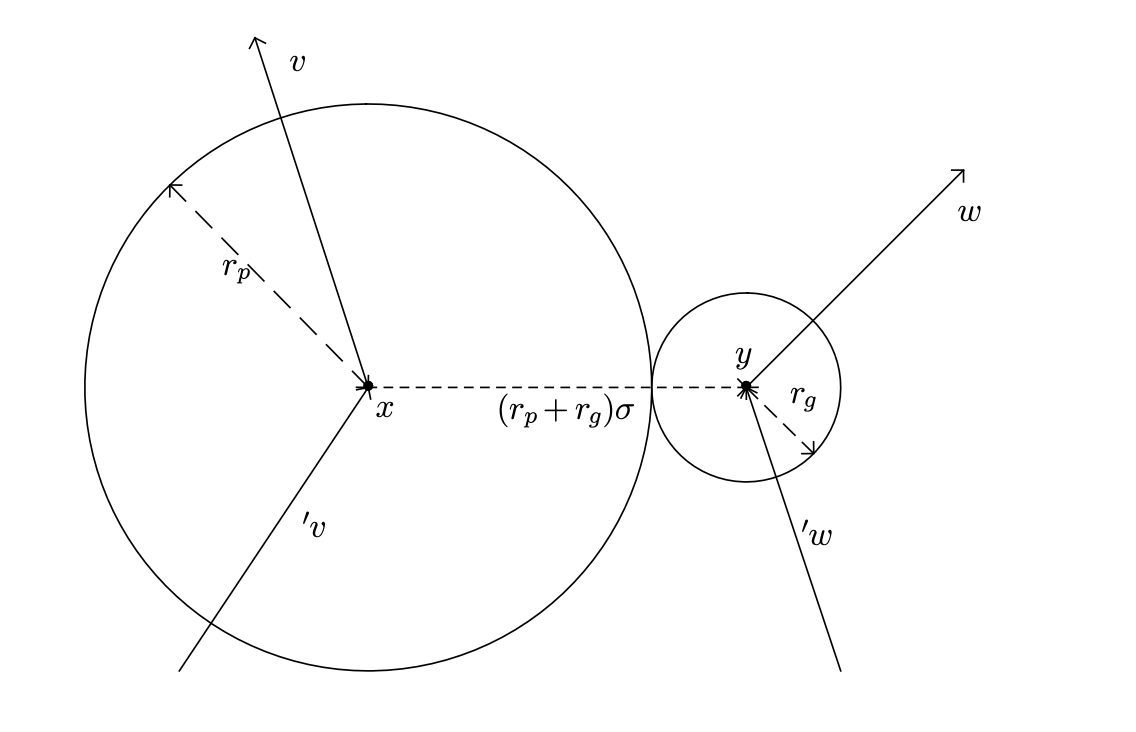}
    \caption{Schematic diagram of collision between different species}
    \label{fig:particle-gas-collision}
\end{figure}

We then introduce the following table of physical parameters (Table \ref{tab:physical parameters}). 
\begin{table}[ht]
\centering
\begin{tabular}{|l|l|}
\hline
Parameter & Meaning \\
\hline
$L$ & size of the container (periodic box) \\
\hline
$T$ & time scale\\
\hline
$N_p$ & number of particles$/L^3$  \\
\hline
$N_g$ & number of gas molecules $/L^3$   \\
\hline
$V_p$ & thermal speed of particles \\
\hline
$V_g$ & thermal speed of gas molecules \\
\hline
$S_p$ & the average collision cross section between particles\\
\hline
$S_g$ & the average collision cross section between gas molecules\\
\hline
$r_p$ & the radius of particles \\
\hline
$r_g$ & the radius of gas molecules \\
\hline
$m_p$ & the mass of each particle \\
\hline
$m_g$ & the mass of each gas molecule \\
\hline
$\eta:=m_g/m_p$ & mass ratio ($\le1$)\\
\hline
$\epsilon:=V_p/V_g$ & thermal speed ratio ($\le1$)\\
\hline
\end{tabular}
\caption{Physical parameters}
\label{tab:physical parameters}
\end{table}
Throughout, subscript $g,p$ refer to gas molecules and particles, respectively. The table of parameters is the same as in \cite{DBGR2017,DBGR2018,DGR2019}. We begin with equation \eqref{Original coupled Enskog-Boltzmann equations}. The dimensionless quantities are given by
\[
\hx=\frac{x}{L},\quad \hat{t}=\frac{t}{T},\quad \hv=\frac{v}{V_p},\quad \hw=\frac{w}{V_g}
\]
and the dimensionless distribution functions are
\[
\hF(\hat{t},\hx,\hv)=\frac{F(t,x,v)}{N_p/V_p^3};\quad \hf(\hat{t},\hx,\hw)=\frac{f(t,x,w)}{N_g/V_g^3}
\]
The collision kernels in the Boltzmann collision integrals are rescaled by

\[
\hat{b}_{g}(\hat{z},\omega)=\frac{b_{g}(z,\omega)}{S_{g}V_g},\qquad \hat{b}_{p}(\hat{z},\omega)=\frac{b_{p}(z,\omega)}{S_{p}V_p}.
\]
The scaling comes from the dimension of the collision integral should match with the left hand side of the Boltzmann equation, see formula 2.2 in \cite{Landau10}.

Summarizing these definitions, we can express the dimensionless collision integrals as follows:
\begin{align*}
\E_0[f,f]&=\frac{N_g^2}{V_g^6}S_{g}V_g^4\hat{\mathcal{E}_0}[\hf,\hf],\\
\E_1[f,F]&=L^2\frac{N_p}{V_p^3}\frac{N_g}{V_g^3}V_gV_p^3\hat{\E_1}[\hf,\hF],\\
\E_2[F,f]&=L^2\frac{N_p}{V_p^3}\frac{N_g}{V_g^3}V_g^4\hat{\E_2}[\hF,\hf],\\
\E_3[F,F]&=\frac{N_p^2}{V_p^6}S_{p}V_p^4\hat{\E_3}[\hF,\hF],
\end{align*}

where dimensionless collision integrals are given by 
\begin{equation*}
\hat{\E_0}[\hf,\hf](\hat{t},\hat{x},\hat{w})=\int_{\RS} [\hf (\hat{x} , \hat{{w}^o}) \hf (\hx-\frac{2r_g}{L}, \hat{w_1}^o) - \hf (\hat{x}, \hw) \hf (\hx+\frac{2r_g}{L}, \hat{w_1})]\hat{b}_{g}( \hw - \hat{w_1}, \n) \dd  \n \dd  \hat{w}_1,
\end{equation*}

\begin{multline*}
\hat{\E_1}[\hf,\hF](\hat{t},\hat{x},\hat{w})=S_g\left(\frac{r_g+r_p}{L}\right)^2\int_{\RS} [\hF (\hat{x} + \frac{r_g+r_p}{L} \n, \hat{v}') \hf (\hx, \hat{w}') - \hF (\hat{x} - \frac{r_g+r_p}{L} \n, \hv) \hf (\hx, \hw)]\\
\times( \epsilon\hv - \hw) \cdot \n H(( \epsilon\hv - \hw) \cdot \n) \dd  \n \dd  \hv
\end{multline*}

\begin{multline*}
\hat{\E_2}[\hF,\hf](\hat{t},\hat{x},\hat{v})=\left(\frac{r_g+r_p}{L}\right)^2\int_{\RS} [\hf (\hat{x} - \frac{r_g+r_p}{L} \n, \hat{w}') \hF (\hx, \hat{v}') - \hf (\hat{x} + \frac{r_g+r_p}{L} \n, \hat{w}) \hF (\hx, \hv)]\\
\times( \epsilon\hv - \hw) \cdot \n H(( \epsilon\hv - \hw) \cdot \n) \dd  \n \dd  \hw
\end{multline*}
\begin{equation*}
\hat{\E_3}[\hF,\hF](\hat{t},\hat{x},\hat{v})=S_p\int_{\RS} [\hF (\hat{x} , \hat{v}^o) \hF (\hx-\frac{2r_p}{L}\n, \hat{v_1}^o) - \hF (\hat{x}, \hv) \hF (\hx+\frac{2r_p}{L}\n, \hat{v_1})]\hat{b}_{p}( \hv - \hat{v_1}, \n) \dd  \n \dd  \hat{v}_1,
\end{equation*}
and the rescaled pre-collision velocities between same species are given by  
\[
\left\{
\begin{aligned}
    \hat{w}^o(\hat{w},\hat{w_1},\n) &:= \widehat{{}^ow}=\hat{w} - \left( \hat{w} - \hat{w_1} \right) \cdot \n\, \n, \\
    \hat{w_1}^o(\hat{w},\hat{w_1},\n) &:=\widehat{{}^ow_1}= \hat{w_1} - \left( \hat{w_1} - \hat{w} \right) \cdot \n\, \n.
\end{aligned}
\right.
\]
 The rescaled pre-collision velocities between different species are given by
\[
\left\{
\begin{array}{l}
     \hat{v}'(\hat{v}, \hat{w}, \n)=\widehat{{}'v}  = \hat{v }- \frac{2 \eta}{1 + \eta} \left( \hat{v} -\frac{1}{\epsilon} \hat {w}
     \right) \cdot \n \n,\\
     \hat{w}'(\hat{v}, \hat{w}, \n)=\widehat{{}'w} = \hat{w} - \frac{2}{1 + \eta} (\hat{w} -  \epsilon\hat{v}) \cdot \n \n,
   \end{array}
\right.
\]
where we recall that $\epsilon:={V_p}/{V_g}$ is the thermal speed ratio, and $\n$ denotes the unit vector along the line connecting the centers of the colliding particles. The dimensionless form of the left-hand side equations should be
\[
\left\{
\begin{aligned}
(\p_t+w\cdot \nabla_x )f(t,x,w)&=\frac{N_g}{V_g^3}\left(\frac{1}{T}\p_{\hat{t}}+\frac{V_g}{L}\hw\cdot \nabla_{\hx}\right)\hf(\hat{t},\hat{x},\hat{w}),\\
(\p_t+v\cdot \nabla_x )F(t,x,v)&=\frac{N_p}{V_p^3}\left(\frac{1}{T}\p_{\hat{t}}+\frac{V_p}{L}\hv\cdot \nabla_{\hx}\right)\hF(\hat{t},\hat{x},\hat{v}).
\end{aligned}\right.
\]
Combining these two parts and dropping the hats on each variable and function, we can write the dimensionless kinetic equations as follows
\begin{equation*}\left\{
\begin{aligned}
\p_t f+\frac{V_gT}{L}w\cdot \nabla_xf&=L^2N_pV_gT \E_1[f,F]+N_gS_gTV_g\E_0[f,f]\\
    \p_t F+\frac{V_pT}{L}v\cdot\nabla_xF&=L^2N_gV_gT\E_2[F,f]+N_pS_pTV_p\E_3[F,F],
    \end{aligned}\right.
\end{equation*}
To simplify the equations, we assume the following conditions
\begin{enumerate}
    \item[(A1)] The velocity scale of gas molecules is the same as that of the particles: $\epsilon= 1$.
    \item[(A2)] We choose the time scale as: $T= L/V_g$.
    \item[(A3)] The length scale satisfies: $L^3N_p= 1$.
    \item[(A4)] The mass ratio is equal to the number ratio, and they are sufficient small: $\frac{N_p}{N_g}=\eta\ll 1$.
    \item[(A5)] We assume the averaged gas molecule cross section $S_g$ such that: ${1}/{\delta}:=N_gLS_g\gg1$.
    \item[(A6)] The collisions between particles can be ignored: $N_pLS_p\ll1$. We will neglect this term in the following discussion.
    \item[(A7)]The radius of gas molecules is sufficient small compared to the length scale: $r_g/L\ll1$. This assumption degenerates Enskog-Boltzmann collision integral $\E_0$ into a classical Boltzmann collision integral $\B$.
\end{enumerate}
In fact, these assumptions are nearly identical to those in \cite{DGR2019}, which are specified for the case of a hard-sphere collision kernel between two species.

Therefore, we arrive at the rescaled kinetic equations:
\begin{equation}\label{dimensionless kinetic equations}
\left\{
\begin{aligned}
    \p_tf+w\cdot\nabla_x f&=\E_1^\eta[f,F]+\frac{1}{\delta}\B[f,f]\\
    \partial_tF + v \cdot \nabla_x F &=\frac{1}{\eta}\E_2^\eta[F,f],
\end{aligned}\right.
\end{equation}
 
where the rescaled Enskog-Boltzmann collision integrals are given by
\begin{align}
    \label{nondimensional Boltzmann collision integral for gas}
    \B[f,f](x,w)=\int_{\RS}  [f (x , w^o) f(x,w_1^o)-f (x , w) f(x,w_1)](w-w_1)\cdot \n
   H(( w-w_1) \cdot \n) \dd  \n \dd  w_1,
\end{align}
\begin{align}
    \label{nondimensional Enskog collision integral for gas}
    &\E_1^\eta[f,F](x,w)=\int_{\RS} a^2 [f (x , w') F (x+a \n, v') - f (x , w) F (x-a \n, v)]
   ( v - w) \cdot \n H(( v - w) \cdot \n)\dd  \n \dd  v,\\
   \label{nondimensional Enskog collision integral for particle}
    &\E_2^\eta[F,f](x,v)=\int_{\RS} a^2 [f (x - a \n, w') F (x, v') - f (x + a \n, w) F (x, v)]
   ( v - w) \cdot \n H(( v - w) \cdot \n) \dd  \n \dd  w,
\end{align}
where the dimensionless radius is defined as
$$a:=\frac{r_p}{L},$$ 
and the pre-collision velocities between different species with mass ratio and velocity scaling are given by 
\begin{equation}\label{scaled collision law between different species}
\left\{ \begin{array}{l}
     v'\equiv v' (v, w, \n,\eta) = v - \frac{2 \eta}{1 + \eta} \left( v -  w
     \right) \cdot \n \n,\\
     w'\equiv w' (v, w, \n,\eta) = w - \frac{2}{1 + \eta} (w -  v) \cdot \n \n,
   \end{array} \right. \end{equation}
and the pre-collision velocities between molecules of the same species with velocity scaling given by
\[ \left\{ \begin{array}{l}
     v^o\equiv v^o (v, v_1, \n) = v - \left( v -  v_1\right) \cdot \n \n\\
     v_1^o\equiv v_1^o(v, v_1, \n) = v_1 -  (v_1-  v) \cdot \n \n .  \end{array} \right. \]

It's easy to check the following collision properties.
\begin{itemize}
    \item Pre-collision velocities are even function with respect to $\n$
    \begin{equation}\label{change the sign of n}
     v' (v,w,\n,\eta) = v'(v,w,-\n,\eta) ; w' (v,w,\n,\eta) = w' (v,w,-\n,\eta).
    \end{equation}
    \item The collision process is reversible: Let $\Phi^\eta(v,w):=(v'(v,w,\n,\eta),w'(v,w,\n,\eta))$, then
    \begin{equation}\label{revisible collision}
         \Phi^\eta\circ \Phi^\eta (v,w)=(v,w).
    \end{equation}
    \item The result above gives the absolute value the determinant of the Jacobian of the change of variable is $1$
    \begin{equation}\label{jacobian of pre-to-post transformation}
    \left|\det \frac{\partial \Phi^\eta(v,w)}{\partial (v, w)} \right|= 1.
    \end{equation}
    \item The collision is specular reflection with respect to norm vector $\n$
    \begin{equation}\label{specular reflection}
     ( v' (v,w,\n,\eta) - w' (v,w,\n,\eta)) \cdot \n = - ( v - w) \cdot \n.
    \end{equation}
\end{itemize}
By simply taking $\eta=1$, we obtain the same properties for the pre-collision velocities $(v^o,v_1^o)$ and $(w^o,w_1^o)$.

At the end of this subsection, we remark that if we adopt a different scaling, for example, if we introduce distinct thermal velocities $V_g$ and $V_p$--- it becomes possible to derive the Vlasov-Stokes, or Vlasov-Navier-Stokes system; see \cite{DBGR2017,DBGR2018}.

\section{Main result: Derivation of the Euler-Vlasov system for thick sprays}\label{sec:Main result}

   Each term in the macroscopic equation can be interpreted from a mesoscopic point of view. By retaining leading-order contributions and carefully estimating higher-order terms, the derived macroscopic equation provides a consistent approximation of the underlying kinetic model while revealing the original of each term from kinetic level.
    
In the following sections, we aim to prove the following theorem, which derive the experimental kinetic-fluid equations \eqref{Vlasov-Euler for thick sprays} from the more fundamental equations—namely, the coupled Enskog-Boltzmann equations, at a formal asymptotic level:

\begin{theorem}\label{thm:derive the thick sprays model}
   Let $f^{\eta,\delta}=f^{\eta,\delta}(t,x,w)$, and $F^{\eta,\delta}=F^{\eta,\delta}(t,x,v)$ solve the following coupled Enskog-Boltzmann equations
   \begin{equation}
       \left\{
      \begin{aligned}
    \p_tf^{\eta,\delta}+w\cdot\nabla_x f^{\eta,\delta}&=\E_1^{\eta}[f^{\eta,\delta},F^{\eta,\delta}]+\frac{1}{\delta}\B[f^{\eta,\delta},f^{\eta,\delta}],\\
    \partial_tF^{\eta,\delta} + v \cdot \nabla_x F^{\eta,\delta} &=\frac{1}{\eta}\E_2^{\eta}[F^{\eta,\delta},f^{\eta,\delta}].
    \end{aligned}\right.
   \end{equation}
   Assume that, for each $(\eta,\delta)$, $f^{\eta,\delta}(w),f^{\eta,\delta}(w)w,f^{\eta,\delta}(w)|w|^2\in L^1(\R_w)$, and $F^{\eta,\delta}\in W^{2,\infty}(\R_x)\cap W^{1,1}(\R_v)$. Moreover, suppose $\p_t f^{\eta,\delta}$, $w\cdot\nabla_xf^{\eta,\delta}$, and $\E_1^\eta[f^{\eta,\delta},F^{\eta,\delta}]$ has uniform bound with respect to $\delta$ and $\eta$. 
   Assume further that there exists $\alpha,\num,u,\theta\in C^1(\mathbb{R_+}\times\R)$, $F\in C^1(\mathbb{R_+},W^{2,\infty}( \R_x,W^{1,1}(\R_v)))$, $f^{\eta,\delta}\to \M\equiv\M[\alpha(t,x)\num(t,x),u(t,x),\theta(t,x)/m_g](w)$, $F^{\eta,\delta}\to F(t,x,v)$, and $\nabla_xF^{\eta,\delta}\to \nabla_xF(t,x,v)$ almost everywhere as $\delta,\eta\to 0$ where $\M$ is Maxwellian distribution
   \begin{equation}\label{maxwellian}
    \M[\alpha\num, u,\theta/m_g](t,x,w):=\frac{\alpha(t,x)\num(t,x)}{(2\pi \theta(t,x)/m_g)^{3/2}}\exp\left(-\frac{m_g|w-u(t,x)|^2}{2\theta(t,x)}\right).
\end{equation}
   Then $(\alpha, \num,u,\theta,F)$ solve the following equations in the weak sense

   \begin{align}
    &\partial_t(\alpha \rho) + \Div_x (\alpha \rho u)=0,\label{gas phase mass equation for thick sprays}\\
    &\p_t (\alpha\rho u)+\Div_x( \alpha\rho u^{\otimes 2})+ \alpha\nabla_x p  =m_g\int_{\R}F(v) D(v-u)\dd v-\Div_x\left(\alpha\rho  \int_{\R}F(v)Q(v-u) \dd v\right)+\mathscr{R},\label{gas phase momentum equation for thick sprays}\\
    &\p_t \left(\alpha\rho E\right)+\Div_x\left( \alpha\rho uE\right)+\Div_x\left( \alpha p u\right)+p \p_t\alpha\label{gas phase energy equation for thick sprays}\\
    &\qquad\qquad =m_g\int_{\R}F(v) D(v-u) v\dd v  -\Div_x\left(\alpha\rho  \int_{\R}F(v)Q(v-u)v \dd v\right) +\mathscr{P},\nonumber\\
    &\p_tF+v\cdot \nabla_xF+\Div_v \left[F\left(D(v-u)+ \frac{4\pi}{3}a^3\nabla_xp\right)\right]=\mathscr{Q},\label{particle equation of thick sprays}\\
    &p=\num\theta,\qquad E=\frac{|u|^2}{2}+\frac{3\theta}{2m_g}, \qquad \rho=m_g n\\
    &\alpha=1-\frac{4\pi}{3}a^3\int_{\R}F(v) \dd v.\label{the volume fraction function}
   \end{align}
   The friction force $D(v-u)\equiv D[\alpha\num,u,\theta](v-u)$ is given by
   \begin{equation}\label{definition of drag force}
   D[\alpha\num,u,\theta](v-u):=  \pi a^2 \alpha\num\frac{\theta}{m_g}\q+\Div_x\left[\alpha\num Q(v-u)\right],
   \end{equation}
   where $q\in \R$ is defined as
  \begin{equation}\label{definition of q}
q(\xi):=\frac{1}{(2\pi)^{3/2}} \int_{\R}\left(\xi-y\right)\left|\xi-y\right|e^{-|y|^2/2}\dd y,
\end{equation}
 $Q$ is a $2$-tensor defined as
\begin{equation}\label{definition of Q}
Q(\xi):=\frac{4\pi}{15}a^3[2\xi^{\otimes 2}+|\xi|^2\id].
\end{equation}
And the remainder terms $\mathscr{P},\mathscr{Q}\in \mathbb{R}$ and $\mathscr{R}\in \R$ are of order $o(a^3)$,  explicitly given by
\begin{equation}\label{definition of remainder terms}
    \begin{aligned}
        \mathscr{P}(t,x):=&2a^2m_g\int_{\R\times\RS} \M[\alpha\num,u,\theta](x,w) [F (x,v)-F (x-a\n,v)+a\n\cdot\nabla_xF (x,v)] 
        \big(( v -w) \cdot \n\big)^2v\cdot \n \\& H(( v -w) \cdot \n)  \dd  v \dd  \n \dd  w -\big((1-\alpha)\num\theta\p_t\alpha\big)(t,x) \\
        \mathscr{Q}(t,x):=&a^2\int_{\RS} F (x, v) \bigg[\M[\alpha\num,u,\theta] (x - a \n, w+2(v-w)\cdot \n \n)-\M[\alpha\num,u,\theta](x,w+2(v-w)\cdot \n \n)\\
        &+a\n\cdot\nabla_x\M[\alpha\num,u,\theta](x,w+2(v-w)\cdot \n \n) - \M[\alpha\num,u,\theta](x + a \n, w)+\M[\alpha\num,u,\theta](x,w)\\
        &+a\n\cdot\nabla_x\M[\alpha\num,u,\theta](x,w) \bigg]( v - w) \cdot \n H(( v - w) \cdot \n) \dd w\dd \n\\
        &+\frac{4\pi}{3m_g}a^3\num\theta\nabla_vF(x,v)\cdot\nabla_x\alpha(t,x)-\frac{4\pi}{3m_g}a^3(1-\alpha)\nabla_vF(x,v)\cdot\big(\nabla_x(\num\theta)\big)(t,x)\\
        \mathscr{R}(t,x):=&2a^2m_g\int_{\R\times\RS}\M[\alpha\num,u,\theta](x,w)[F(x-a\n,v)-F(x,v)+a\n\cdot\nabla_xF(x,v)][(v-w)\cdot \n]^2\n\\
        &H((v-w)\cdot \n)\dd v\dd w\dd \n+\big((1-\alpha)\num\theta\nabla_x\alpha\big)(t,x),
    \end{aligned}
\end{equation}
   \end{theorem}
  The theorem is primarily based on the following two propositions, which are proved in Section \ref{sec:Particle equation} and Section \ref{sec:Gas phase equations}, respectively. Proposition \ref{prop:weak formula for drag force} aims to derive the Vlasov-type equation \eqref{particle equation of thick sprays} from \eqref{dimensionless kinetic equations}.
\begin{proposition}\label{prop:weak formula for drag force} For any test function $\varphi(v)\in C^\infty(\R)$, and Enskog-Boltzmann collision integral $\mathcal{E}_2^\eta[F,f]$ defined in \eqref{nondimensional Enskog collision integral for particle}, $f^{\eta,\delta},F^{\eta,\delta}$ satisfy the assumptions shown in Theorem \ref{thm:derive the thick sprays model}, then we have the following weak formulation:
\begin{equation}\label{weak formulation of E_2}
\lim_{\eta,\delta\to0}\frac{1}{\eta}\int_{\R}\E_2^\eta[F^{\eta,\delta},f^{\eta,\delta}]\varphi(v)\dd v= \int_{\R}\left[F(v)\left(D(v-u)+ \frac{4\pi}{3}a^3\nabla_xp\right)\right]\cdot\nabla_v\varphi(v)\dd v  +\int_{\R}\mathscr{Q}\varphi(v)\dd v
\end{equation}
where the drag force $\Gamma$ is defined in \eqref{definition of drag force}, and remain term $\mathscr{Q}$ is defined in \eqref{definition of remainder terms}.
\end{proposition}
With the help of this proposition, we multiply the test function $\varphi(v)$ on both sides of Enskog-Boltzmann equation for particles in \eqref{dimensionless kinetic equations}, and integrate it with respect to $v\in\R$. As $\eta\to 0$, the distribution function $F$ should satisfies
\[
\int_{\R}\p_t\varphi(v) F(v)\dd v+\int_{\R}\nabla_x\varphi(v)\cdot v F(v)\dd v+\int_{\R}\nabla_v\varphi(v)\cdot\left[F(v)\left(D(v-u)+ \frac{4\pi}{3}a^3\nabla_xp\right)\right]\dd v=-\int_{\R}\mathscr{Q}\varphi(v)\dd v,
\]
 which corresponds to a Vlasov-type equation in the weak sense. Proposition \ref{prop:gas phase equations} aims to obtain the gas phase equations \eqref{gas phase mass equation for thick sprays},\eqref{gas phase momentum equation for thick sprays}, and \eqref{gas phase energy equation for thick sprays}---Euler-like equations---from coupled Enskog-Boltzmann equations \eqref{dimensionless kinetic equations}.
\begin{proposition}\label{prop:gas phase equations}
Let $f^{\eta,\delta},F^{\eta,\delta}$ satisfy the assumptions shown in Theorem \ref{thm:derive the thick sprays model}, then the parameter functions $\alpha,\num,u,\theta$ of $\M[\alpha\num, u,\theta/m_g]$, the limit function of $f^{\eta,\delta}$, obeys the following equations:
\begin{equation}
    \left\{
    \begin{aligned}
    &\partial_t(\alpha \rho) + \Div_x (\alpha \rho u) =0\\
    &\p_t (\alpha\rho u)+\nabla_x( \alpha\rho u^{\otimes 2})+ \alpha\nabla_x(\num\theta )=m_g\int_{\R}  D(v-u)F(v)\dd  v-\Div_x\left(\alpha\rho\int_{\R}Q(v-u) F(v)\dd v\right) +\mathscr{R},\\
     &\p_t \left(\alpha \left(\frac{\rho|u|^2}{2}+\frac{3\num\theta}{2} \right)\right)+\Div_x\left( \alpha u\left(\frac{\rho|u|^2}{2}+\frac{3\num\theta}{2} \right)\right)+\Div_x(\alpha \num \theta u)+\num\theta\p_t\alpha\\
     &\qquad=m_g\int_{\R}F(v) D(v-u)\cdot v \dd  v  -\Div_x \left(\alpha\rho\int_{\R}  Q(v-u) vF(v)  \dd  v\right) +\mathscr{P},
    \end{aligned}\right.
\end{equation}
where the friction force $D(v-u)$ is defined in \eqref{definition of drag force}, and $Q(v-u) $ is defined in \eqref{definition of Q}. $rho$ is local mass density of gas defined by $\rho(t,x)=m_g \num(t,x)$.
\end{proposition}

This theorem explains why the volume fraction $\alpha$ appears outside the pressure gradient term $\nabla_xp$: the microscopic drag force inherently contains a pressure contribution. When this term interacts with density distribution of the particles, it yields an additional volume fraction function, specifically, $\frac{4\pi}{3}a^3\int_{\R}F(v)\dd v\nabla_xp$. If we separate this volume effect and friction force, and combine it with the pressure gradient $\nabla_x(\alpha p)$ that arises from the free transport (i.e., the gas phase's self-interaction), we are naturally led to a term of the form $\alpha\nabla_x(\alpha p)$.
In fact, under the assumption that $\alpha$  is close to $1$, it can be approximated by $\alpha\nabla_x p+o(a^3)$, which justifies writing the effective pressure gradient as $\alpha\nabla_xp$ at leading order. This formal argument reveals the microscopic origin of the pressure modification observed in thick sprays models.  

The appearance of the coefficient $\frac{4}{3}\pi a^3$ is not coincidental. In fact, if this problem is treated in dimension $d$, the corresponding coefficient would precisely be the volume of a $d$-dimensional ball. See Appendix \ref{apd:general dimension}.

 The leading-order term $q$ of the friction force $D$ can be expressed as 
 \[
 \q=\bar{q}\left(\frac{|v-u|}{\sqrt{\theta/m_g}}\right)(v-u),
 \]
 where the explicit formula for $\bar{q}$ can be found in Lemma 3.3 of \cite{DGR2019}. However, the delocalization effect introduce not only a modification to the pressure term but also an additional friction term $\Div_x(\alpha\rho Q(v-u))$. If we write friction force $D$ in terms of a viscous tensor $\bar{D}(v-u)$ acting on the relative velocity $v-u$, that is,
 \[
 D[\alpha\num,u,\theta]=\bar{D}[\alpha\num,u,\theta](v-u),
 \]
 we observe that the viscous tensor is no longer a scalar. Instead, the resulting friction force may deviate toward the direction of  $\nabla_x(\alpha \num)$, and its structure also depends on the compressibility of the gas phase.
 The explicit expression for viscous tensor $\bar{D}(v-u)$ can be found in Appendix \ref{apd:more on the friction force}. Consequently, the resulting friction force may no longer strictly align with the direction with $v-u$. 
 Nonetheless, since this deviation $\Div_x(\alpha\rho Q(v-u))$ is arises at higher order $O(a^3)$ compared to the leading-order friction term $q$, it may be regraded as a small perturbation.

Another way to derive the gas phase equation is combining the following proposition with proposition \ref{prop:weak formula for drag force}. 
\begin{proposition}\label{prop:weak formula for enskog collision integral}
For Enskog-Boltzmann collision integrals defined in \eqref{nondimensional Enskog collision integral for gas},\eqref{nondimensional Enskog collision integral for particle} obeys the collision laws \eqref{scaled collision law between different species}. For any test function $\varphi(x,v)\in C^\infty(\R\times\R)$, we have the following identity

\begin{multline}\label{weak formula for summation of collision integrals}
\eta\int_{\R} \E_1^\eta[f,F](w) \varphi(x,w)\dd w+\int_{\R}\E_2^\eta[F,f](v)\varphi(x,v)\dd v+\eta \Div_x  I [f,F;\varphi]\\
=\int_{\RS\times\R}a^2 f(x,w)F(x-a\n,v)(v-w)\cdot \n H((v-w)\cdot \n)\bigg[\varphi(x-a\n,w')+\varphi(x,v')\\
-\varphi(x-a\n,w)-\varphi(x,v)\bigg]\dd \n \dd v\dd w,
\end{multline}
where 
\begin{multline*}
     I [f,F;\varphi]:=\int_{\RS\times\R}\int^a_0a^2 \n(v-w)\cdot \n H((v-w)\cdot \n)(\varphi(x+s \n,w')-\varphi(x+s\n,w)\big)\\f (x+s \n , w) F (x-(a-s) \n, v)
    \dd \n \dd v\dd w\dd s .
\end{multline*}
\end{proposition}
The proof of the Proposition \ref{prop:weak formula for enskog collision integral} can be found in Appendix \ref{apd:proof of prop 3}.
In fact, the terms $\alpha\rho\int_{\R}FQ \dd v$ and $\alpha\rho\int_{\R}FQ v\dd v$ correspond precisely to the leading-order contributions (of order $O(a^3)$) of $I[f,F;\varphi]$ for $\varphi(x,\xi)=\xi$ and $\varphi(x,\xi)=\frac{|\xi|^2}{2}$ respectively.

The results for general delocalized collision integral in the single-species case can be found in \cite{CCG2024}, and the corresponding result for the Boltzmann collision integral for binary mixture can be found in \cite{DBGR2017,DBGR2018,DGR2019}.

This proposition illustrates more clearly the origin of the pressure imbalance, and explain that the pressure imbalance is inevitable when we introduce ``delocalized collision'' models. Although, these pressures do not balance locally, their integral over $x\in\R$ causes the ``extra pressure'' to vanishes. Moreover, in non-equilibrium gas-particle systems, this phenomena occurs, see, for example \cite{GHS2004,DM2010}. It's referred to as interfacial pressure. 


\section{Proof of proposition \ref{prop:weak formula for drag force}}\label{sec:Particle equation}

In this section, we start from equation  \eqref{dimensionless kinetic equations}, with the goal of deriving equation \eqref{particle equation of thick sprays} which incorporates the drag force \eqref{definition of drag force} and the volume fraction \eqref{the volume fraction function}. We first derive the Maxwellian distribution from the Enskog-Boltzmann equation for gas phase. Next, we obtain the leading order drag force in the following subsection. In the final subsection, we derive the higher-order correction to drag force, which includes contribution from the volume fraction and pressure.

\subsection{The leading order gives the Maxwellian distribution}
We assume that $(f^{\eta,\delta},F^{\eta,\delta})$ solves the equation \eqref{dimensionless kinetic equations}. Hence
\[
\B[\fed,\fed]=\delta(\p_t\fed+w\cdot\nabla_x\fed-\E_2^\eta[\Fed,\fed]).
\]
Since right-hand side is uniformlly bounded in $\delta$, as $\eta,\delta\to 0$, 
the limit function $f$ should satisfy
\begin{equation}\label{f conserves all quantities}
\B[f,f](t,x,w)=0.
\end{equation}
It implies that $f$ is a Maxwellian distribution; see for instance, Chapter 3 of \cite{cercignani1994}. Moreover, we introduce the following macroscopic quantities
\begin{align}
&\rho(t,x)=m_g n(t,x)=\frac{H(\alpha(t,x))}{\alpha(t,x)}\int_{\R}m_g f(t,x,w)\dd w,&\mathrm{the~local~mass~density~of~gas}\label{the local density of gas}\\
&u(t,x)=\frac{H(\alpha(t,x)\num(t,x))}{\alpha(t,x)\num(t,x)}\int_{\R} wf(t,x,w)\dd w,&\mathrm{the~bulk~velocity~of~gas}\label{the bluk velocity of gas}\\
&\theta(t,x)=\frac{H(\alpha(t,x)\num(t,x))}{\alpha(t,x)\num(t,x)}\int_{\R}m_g \frac{|w-u(t,x)|^2}{3}f(t,x,w)\dd w,&\mathrm{the~internal~energy~of ~gas}\label{the internal energy of gas}
\end{align}
then $f^\eta(t,x,w)$ should be of the form of Maxwellian distribution $\M[\alpha\num,u,\theta/m_g](w)$ in \eqref{maxwellian},
where the volume fraction function $\alpha=\alpha(t,x)$ is defined by
\begin{equation}\label{definition of volume fraction}
    \alpha(t,x):=1-\frac{4\pi}{3}a^3\int_{\R}F(t,x,v)\dd v,
\end{equation}
to describe the volume occupation of gas which is close to $1$ but not negligible. 
These quantities are well-defined thanks to the strict  positivity of volume fraction $\alpha$. This can be established in the simplified case (with linear viscous force). We refer the reader to the Proposition 2.1 in \cite{BDD2023}.


\subsection{Drag force from Boltzmann collision integral}
As in \cite{DGR2019}, we can construct the term $q(v-u)F$ from $\frac{1}{\eta} E_2[F, f]$. This motivates us to remove the delocalization effect on $f$ by using the approximation
\[
f^{\eta,\delta}(x\pm a\n, w)=f^{\eta,\delta}(x,w)\pm a\n\cdot \nabla_x f^{\eta,\delta}(x,w)+o(a).
\]
Asymptotically, the Enskog-Boltzmann collision integral can be rewritten as 

\begin{equation}\label{expansion on enskog collision integral}
   \E_2^\eta [F^{\eta,\delta}, f^{\eta,\delta}] = \mathcal{R}_1^\eta [F^{\eta,\delta}, f^{\eta,\delta}] -\mathcal{R}_2^\eta[F^{\eta,\delta},f^{\eta,\delta}]+
   o(a^3),  
\end{equation} 
where $\mathcal{R}_2^\eta[F^{\eta,\delta},f^{\eta,\delta}]$ and $\mathcal{R}_2^\eta[F^{\eta,\delta},f^{\eta,\delta}]$ are defined by
    \begin{align}\label{definition of R}
    \mathcal{R}_1^\eta[F^{\eta,\delta},f^{\eta,\delta}](t,x):=&a^2 \int_{\RS} [f^{\eta,\delta} (x, w') F^{\eta,\delta} (x, v') -
    f^{\eta,\delta} (x, w) F^{\eta,\delta} (x, v)] ( v - w) \cdot \n H(( v -w) \cdot \n)\dd  \n \dd  w ,\\
   \mathcal{R}_2^\eta[F^{\eta,\delta},f^{\eta,\delta}](t,x):=&a^3 \int_{\RS} \n \cdot [\nabla_x f^{\eta,\delta} (x, w') F^{\eta,\delta} (x, v') +
   \nabla_x f^{\eta,\delta} (x, w) F^{\eta,\delta} (x, v)] ( v - w) \cdot \n \nonumber\\
   &\qquad \times H(( v -w) \cdot \n)\dd  \n \dd  w.\nonumber
    \end{align}
Hence, we define 
\begin{equation}\label{remainder Q_1}
\mathscr{Q}_1^{\eta,\delta}:=\E_2^\eta [F^{\eta,\delta}, f^{\eta,\delta}] - \mathcal{R}_1^\eta [F^{\eta,\delta}, f^{\eta,\delta}] +\mathcal{R}_2^\eta[F^{\eta,\delta},f^{\eta,\delta}],
\end{equation}
and $\mathscr{Q}_1:=\lim_{\eta,\delta\to0}\mathscr{Q}_1^{\eta,\delta}$.
The computation of leading term follows the same approach as in \cite{CDS2012}. We summarize only the main idea here; for the inelastic case, one may refer to \cite{CDS2012}, also in \cite{DGR2019}.
For a test function $\varphi(v) \in C^2 (\mathbb{R}^3)$, with $\nabla\varphi,\nabla^2\varphi\in L^\infty(\R)$, we apply the change of variables $(v,w,\n)\mapsto (v',w',-\n)$ in the gain term of  $\mathcal{R}_1^\eta$. Properties \eqref{revisible collision}, \eqref{jacobian of pre-to-post transformation}, and \eqref{specular reflection} guarantee that
\begin{multline*}
   \int_{\R\times\RS}  [ f^{\eta,\delta} (w') F^{\eta,\delta} (v') -  f^{\eta,\delta} (w) F^{\eta,\delta} (v)] ( v -
  w) \cdot \n H(( v -w) \cdot \n) \varphi (v) \dd  v \dd  \n \dd  w  \\
  =\int_{\R\times\RS}  f^{\eta,\delta} (w) F^{\eta,\delta} (v) [\varphi(v')-\varphi(v)]( v -
  w) \cdot \n H(( v -w) \cdot \n)  \dd  v \dd  \n \dd  w. 
\end{multline*}
The fundamental theorem of calculus, and \eqref{scaled collision law between different species} show that
\begin{equation}\label{cmp:taylor expansion on test function}
\varphi(v')-\varphi(v)=(v'-v)\cdot \nabla_v\varphi(v)+o(v'-v)=-\frac{2 \eta}{1 + \eta} ( v - w) \cdot \n \n\cdot\nabla_v \varphi(v)  +o(\eta),
\end{equation}
hence
\begin{multline}\label{weak formulation of B_2} 
\int_{\R}\mathcal{R}_1^\eta[F^{\eta,\delta},f^{\eta,\delta}] (x, v)\varphi(v)\dd v \\
= -\frac{2a^2 \eta }{1 + \eta} \int_{\R}\left(
   F^{\eta,\delta} (x, v) \int_{\mathbb{S}^2}  f^{\eta,\delta} (w) (( v - w) \cdot \n)^2H(( v - w) \cdot \n) \n
   \dd  w \dd  \n \right)\cdot\nabla_v\varphi(v)\dd v  + o (\eta),
\end{multline}
which is the weak formulation for Boltzmann collision integral $\mathcal{R}_1^\eta$. After integrating in $\n$ (see, for instance, section 16.8 of \cite{chapman1970}, or in the Appendix \ref{apd:two integrals}), we obtain
\begin{equation} \label{drag force term before integrating on sphere}
\int_{\mathbb{S}^2}  \n (( v - w) \cdot \n)^2H(( v - w) \cdot \n)f^{\eta,\delta}(w)\dd  w \dd  \n = \frac{\pi}{2} 
    \int  ( v - w)|
    v - w | f^{\eta,\delta}(w) \dd  w.  
     \end{equation}
Since $f^{\eta,\delta}\to f^\eta$ takes the form of a Maxwellian $\M[\alpha\num,u,\theta/m_g]$, 
\begin{equation}
 \int_{\mathbb{R}^3} (v - w) | v - w | \mathcal{M} [\alpha\num,u,\theta/m_g] (w) \dd  w=\alpha\num\frac{\theta}{m_g}\q
\end{equation}
where $q$ is defined in \eqref{definition of q}. The explicit formulas for $q$ can be found in Lemma 3.3 of \cite{DGR2019}, and their asymptotic behavior is discussed in \cite{CMS2025}. 
As a result, we obtain the leading-order friction force as
\begin{equation}\label{drag force term}
    \lim_{\delta,\eta\to 0}\int_{\R}\mathcal{R}_1^\eta[F^{\eta,\delta},f^{\eta,\delta}] (x, v)\varphi(v)\dd v \\
= -\pi a^2\alpha\num\frac{\theta}{m_g} \int_{\R}F(v)\q\cdot\nabla_v\varphi(v)\dd v  ,
\end{equation}
\subsection{Volume fraction from delocalization effect}
We repeat the same procedure as in the previous subsection to derive the weak formulation of $\mathcal{R}_2^\eta$. For any test function $ \varphi \in C^2(\mathbb{R}^3)$ with $\nabla \varphi, \nabla^2 \varphi \in L^\infty(\mathbb{R}^3)$, we perform the change of variables $(v, w, n) \mapsto (v'(v,w,n,\eta), w'(v,w,n,\eta), -n) $ in the gain term (first half part) of $ \mathcal{R}_2^\eta$. Consequently, we obtain:  
\begin{multline*}
    \int_{\R}\mathcal{R}_2^\eta[F^{\eta,\delta},f^{\eta,\delta}] (x, v)\varphi(v)\dd v \\
    = a^3\int_{\R\times\RS} \n\cdot\nabla_x f^{\eta,\delta} (w) F^{\eta,\delta} (v) [-\varphi(v')+\varphi(v)]( v -
  w) \cdot \n H(( v -w) \cdot \n)  \dd  v \dd  \n \dd  w.
\end{multline*}
Then, the approximation \eqref{cmp:taylor expansion on test function} shows 
\begin{multline}
    \int_{\R}\mathcal{R}_2^\eta[F^{\eta,\delta},f^{\eta,\delta}] (x, v)\varphi(v)\dd v \\
    = \frac{2a^3 \eta }{1 + \eta}  \int_{\R}\left(
   F^{\eta,\delta} (x, v)\int_{\mathbb{S}^2} \n \cdot \nabla_x f^{\eta,\delta} (w) ((  v - w) \cdot \n)^2H((v-w)\cdot \n) \n
   \dd  w \dd  \n \right)\cdot\nabla_v\varphi(v)\dd v  + o (\eta).
\end{multline} 
Note the following identity (see, for instance, section 16.8 of \cite{chapman1970}, or in the Appendix \ref{apd:two integrals}) 
\begin{equation}\label{4-tensor integral}
2a^3\int_{|\n|=1,\n\cdot \xi>0}(\n\cdot \xi)^2\n^{\otimes 2}\dd \n=\frac{4\pi}{15}a^3(2\xi^{\otimes 2}+|\xi|^2\mathrm{Id})=:Q(\xi).
\end{equation}
It follows
\begin{multline}\label{cpt:pressure term}
\int_{\RS} \n \cdot \nabla_x f^{\eta,\delta} (w) ((  v - w) \cdot \n)^2H((v-w)\cdot \n) \n
   \dd  w \dd  \n \\
   = \frac{2 \pi}{15} \left(
   \Div_x\int 2 (  v - w)^{\otimes 2} f^{\eta,\delta}(w) \dd  w + \nabla_x \int |
     v - w |^2 f^{\eta,\delta}(w) \dd  w \right) .
\end{multline}

We decompose $v-w$ as $(v-u)+(u-w)$. The component depending on $v-u$ contributes to the pressure term, while the component depending on $u-w$ gives rise to a higher order friction force. In other words,
\begin{multline*}
   \lim_{\eta,\delta\to0}\int_{\R} (  v - w)^{\otimes 2} f^{\eta,\delta}(w) \dd  w =(v-u)^{\otimes 2}\int_{\R}\M[\alpha\num,u,\theta/m_g](w)\dd w+\int_{\R}(u-w)^{\otimes 2}\M[\alpha\num,u,\theta/m_g](w)\dd w\\
   +(v-u)\otimes\int_{\R}(u-w)\M[\alpha\num,u,\theta/m_g](w)\dd w+\int_{\R}(u-w)\M[\alpha\num,u,\theta/m_g](w)\dd w\otimes (v-u).
\end{multline*}
Since, $\M[\alpha\num,u,\theta/m_g](w)$ is a radial function centered at $w=u$, we have
\[
\int_{\R}(w-u)\M[\alpha\num,u,\theta/m_g](w)\dd w=0.
\]
Therefore,
\begin{equation}\label{2-tensor integral}
  \lim_{\delta,\eta\to 0} \int_{\R} (  v - w)^{\otimes 2} f^{\eta,\delta}(w) \dd  w =\alpha\num(v-u)^{\otimes 2}+\alpha\num\frac{\theta}{m_g}\mathrm{Id}.
\end{equation}
Take the trace of the identity \eqref{2-tensor integral} for both sides, we obtain
\begin{equation}\label{square integral}
   \lim_{\eta,\delta\to0}\int_{\R} |v - w|^2 f^{\eta,\delta}(w) \dd  w =\alpha\num|v-u|^{2}+3\alpha\num\frac{\theta}{m_g}.
\end{equation}
Substituting \eqref{2-tensor integral} and \eqref{square integral} into \eqref{cpt:pressure term}, we deduce that
\begin{multline*}
 \lim_{\eta,\delta\to0}\int_{\RS} \n \cdot \nabla_x f^{\eta,\delta} (w) ((  v - w) \cdot \n)^2H((v-w)\cdot \n) \n
   \dd  w \dd  \n \\
   =\frac{2\pi}{3m_g}\nabla_x(\alpha \num \theta)+\frac{2 \pi}{15m_g} \left[\Div_x  \left( 2\alpha\num(v - u)^{\otimes 2}\right) + \nabla_x \left(\alpha \num| v - u |^2\right) \right] .
\end{multline*}
Next, we aim to eliminate the double appearance of $\alpha$ by discarding terms of order $O(a^3)$. In other words, we will show that the leading-order approximation of $\alpha \nabla_x(\alpha \num \theta)$ is $\alpha \nabla_x(\num \theta)$.
Since $1-\alpha=O(a^3)$ and $\nabla_x\alpha=O(a^3)$, we define
\[
\mathscr{Q}_2:=\frac{1}{m_g}\frac{4\pi}{3}a^3\left[\num\theta\nabla_vF\cdot\nabla_x\alpha-(1-\alpha)\nabla_vF\cdot\nabla_x(\num\theta)\right]=o(a^3).
\]
Therefore, by applying the production rule repeatedly, we have
\begin{align*}
    &\frac{4\pi}{3}a^3\int_{\R}F(v)\nabla_x(\num\theta)\cdot\nabla_v\varphi(v)\dd v-\int_{\R}\mathscr{Q}_2(v)\varphi(v)\dd v\\
=&-\frac{4\pi}{3}a^3\int_{\R}\nabla_vF(v)\cdot\nabla_x(\num\theta)\varphi(v)\dd v-\int_{\R}\mathscr{Q}_2(v)\varphi(v)\dd v\\
=&\frac{4\pi}{3}a^3\int_{\R}\nabla_vF(v)\cdot[-\nabla_x(\num\theta)-\num\theta\nabla_x\alpha+(1-\alpha)\nabla_x(\num\theta)]\varphi(v)\dd v\\
=&-\frac{4\pi}{3}a^3\int_{\R}\nabla_vF(v)\cdot\nabla_x(\alpha\num\theta)\varphi(v)\dd v\\
=&\frac{4\pi}{3}a^3\int_{\R}F(v)\nabla_x(\alpha\num\theta)\cdot \nabla_v\varphi(v)\dd v.
\end{align*}

Then the weak formulation of $\mathcal{R}_2[F,f]$ can be reacted as
\begin{multline}\label{pressure term}
\lim_{\eta,\delta\to0}\frac{1}{\eta}\int_{\R}\mathcal{R}_2^\eta[F^{\eta,\delta},f^{\eta,\delta}] (v)\varphi(v)\dd v =  \int_{\R}F(v)\left(\frac{1}{m_g}\frac{4\pi}{3}a^3\nabla_x(\num\theta)+\Div_x[\alpha\num Q(v-u)]\right)\cdot\nabla_v\varphi(v)\dd v\\
    -\int_{\R}\mathscr{Q}_2\varphi(v)\dd v.
\end{multline}

Finally, we conclude the results by combining \eqref{remainder Q_1}, \eqref{drag force term} and \eqref{pressure term}:
\begin{multline*}
\lim_{\eta,\delta\to0}\frac{1}{\eta}\int_{\R}\E_2^\eta[F^{\eta,\delta},f^{\eta,\delta}]\varphi(v)\dd v= -\int_{\R}\left[F(v)\left(D(v-u)+ \frac{1}{m_g}\frac{4\pi}{3}a^3\nabla_x(\num\theta)\right)\right]\cdot\nabla_v\varphi(v)\dd v  \\
+\int_{\R}\mathscr{Q}\varphi(v)\dd v,
\end{multline*}
where the remainder term $\mathscr{Q}$ is defined by $\mathscr{Q}:=\mathscr{Q}_1+\mathscr{Q}_2$.

\section{Proof of proposition \ref{prop:gas phase equations}}\label{sec:Gas phase equations}
In this section, we use the moment method to derive the equations of gas phase by taking moments of kinetic equation of gas in \eqref{dimensionless kinetic equations} with respect to the functions $m_g,m_gw$, and $\frac{m_g|w|^2}{2}$.



\subsection{Gas phase mass equation}
Once we integrate the gas phase equation in \eqref{dimensionless kinetic equations} with respect to $w$, we obtain  
\[
m_g\partial_t \int_{\R} f^{\eta,\delta}\,\dd w + m_g\Div_x \int_{\R} w f^{\eta,\delta}\,\dd w = m_g\int_{\R}\E_1^\eta[f^{\eta,\delta},F^{\eta,\delta}]\dd w+m_g\int_{\R}\mathcal{B}^\eta[f^{\eta,\delta},f^{\eta,\delta}]\dd w=0.
\]  
Let $\eta,\delta\to 0$, we obtain
\[
m_g\partial_t\int_{\R}\M[\alpha\num,u,\theta/m_g](w)\dd w + m_g\Div_x \int_{\R}\M[\alpha\num,u,\theta/m_g](w) w\dd w = 0.
\]
The definitions \eqref{the local density of gas},\eqref{the bluk velocity of gas} and assumption \eqref{maxwellian} yield 
\begin{equation}\label{conservation of mass for gas phase}
\partial_t(\alpha \rho) + \Div_x (\alpha \rho u)=0.
\end{equation}
The last identity holds because the Enskog-Boltzmann collision integral for two species preserves mass for suitable functions $f$ and $F$. The general form of this property can be found in \cite{CCG2024}.

\subsection{Gas phase momentum equations}
Similarly, we multiply gas the phase equation in \eqref{dimensionless kinetic equations} by $m_gw$, and then integrate it with respect to $w$, then the left-hand side should be
\begin{multline}
\lim_{{\eta,\delta}\to 0}m_g\int_{\R}w\p_t f^{\eta,\delta}(w)\dd w+m_g\int_{\R}w^{\otimes 2}:\nabla_x f^{\eta,\delta}(w)\dd w\\
=\p_t\int_{\R}m_g\M[\alpha\num,u,\theta/m_g](w)w\dd w+\Div_x\int_{\R}m_g\M[\alpha\num,u,\theta/m_g](w)w^{\otimes 2}\dd w.
\end{multline}
We split $w^{\otimes 2}$ into  $w-u$ and $u$. Noting that $\M[\alpha\num,u,\theta/m_g](w)(w-u)$ is an odd function with respect to $w-u$,
\[
\int_{\R}m_g\M[\alpha\num,u,\theta/m_g](w)w^{\otimes 2}\dd w=u^{\otimes 2}\int_{\R}m_g\M[\alpha\num,u,\theta/m_g](w)\dd w+\int_{\R}m_g\M[\alpha\num,u,\theta/m_g](w)(w-u)^{\otimes 2}\dd w,
\]
The limit becomes
\begin{equation}\label{transport part of momentum equation of the gas phase}
    \lim_{{\eta,\delta}\to 0}\int_{\R}m_gw\p_t f^{\eta,\delta}(w)\dd w+\int_{\R}m_gw^{\otimes 2}:\Div_x f^{\eta,\delta}(w)\dd w\\
=\p_t (\alpha\rho u)+\nabla_x( \alpha\rho u^{\otimes 2})+ \nabla_x(\alpha\num\theta )
\end{equation}

We proceed with a similar computation shown in Section \ref{sec:Particle equation}. The property \eqref{revisible collision} allows us combining the gain term and lose term in collision integral $\E_1[f,F]$:
\begin{multline*}
    \int_{\R}\E_1^\eta[f^{\eta,\delta},F^{\eta,\delta}]w\dd w  =\int_{\R\times\RS} a^2 f^{\eta,\delta} (x,w) F^{\eta,\delta} (x-a\n,v) (w'-w)( v -
  w) \cdot \n H(( v -w) \cdot \n)  \dd  v \dd  \n \dd  w\\
  =\frac{2a^2}{1+\eta}\int_{\R\times\RS}  f^{\eta,\delta} (x,w) F^{\eta,\delta} (x-a\n,v) \big(( v -
  w)\cdot \n\big)^2\n H(( v -w) \cdot \n)  \dd  v \dd  \n \dd  w.
\end{multline*}
By approximating $F^{\eta,\delta}(x-a \n,\xi)$ as $F^{\eta,\delta}(x,\xi)-a\n\cdot \nabla_xF^{\eta,\delta}(x,\xi)+o(a)$ in $\int_{\R}\E_1[f^{\eta,\delta},F^{\eta,\delta}]w\dd w$, one has
\begin{multline*}
    \int_{\R}\E_1^\eta[f^{\eta,\delta},F^{\eta,\delta}]w\dd w  =\frac{2a^2}{1+\eta}\int_{\R\times\RS} f^{\eta,\delta} (w) F^{\eta,\delta} (v) \big(( v -
  w) \cdot \n\big)^2\n H(( v -w) \cdot \n)  \dd  v \dd  \n \dd  w\\
  -\frac{2a^3}{1+\eta}\int_{\R\times\RS}  f^{\eta,\delta} (w)  \big(( v -
  w)\cdot \n\big)^2\n^{\otimes 2} \nabla_x F^{\eta,\delta} (v) H(( v -w) \cdot \n)  \dd  v \dd  \n \dd  w+\mathscr{R}_1^{\eta,\delta},
\end{multline*}
where $\mathscr{R}_1^{\eta,\delta}=o(a^3)$ is given by
\begin{multline*}
    \mathscr{R}_1^{\eta,\delta}:=\frac{2a^2}{1+\eta}\int_{\R\times\RS}  f^{\eta,\delta} (x,w) [F^{\eta,\delta} (x-a\n,v)-F^{\eta,\delta} (x,v)+a\n\cdot\nabla_xF^{\eta,\delta} (x,v)] \big(( v -
  w)\cdot \n\big)^2\n \\
  H(( v -w) \cdot \n)  \dd  v \dd  \n \dd  w.
\end{multline*}

Then, integrating with respect to $\n\in\mathbb{S}^2$, and using identities \eqref{drag force term before integrating on sphere}, and \eqref{4-tensor integral}, we obtain
\begin{multline*}
    \int_{\R}\E_1^\eta[f^{\eta,\delta},F^{\eta,\delta}]w\dd w  =\frac{\pi a^2}{1+\eta}\int_{\R\times\R} F^{\eta,\delta} (v)f^{\eta,\delta}(w)|v-w|(v-w)  \dd  v  \dd  w\\
  -\frac{4\pi a^3}{15(1+\eta)}\int_{\R\times\R} [2(v-w)^{\otimes 2}+|v-w|^2\id] \nabla_xF^{\eta,\delta}(v) f^{\eta,\delta} (w)   \dd  v \dd  w+\mathscr{R}_1^{\eta,\delta}.
\end{multline*}

After that, after taking ${\eta,\delta}\to 0$, identities \eqref{drag force term}, \eqref{2-tensor integral}, and \eqref{square integral} yield
\begin{multline*}
   \lim_{{\eta,\delta}\to 0} \int_{\R}\E_1^\eta[f^{\eta,\delta},F^{\eta,\delta}]m_gw\dd w  =\pi a^2\alpha\num\theta\int_{\R} F (v) \q  \dd  v \\
  -\int_{\R} \left\{\frac{4\pi}{3}a^3\alpha \num\theta\id+\alpha\rho Q(v-u)]\right\}\nabla_xF(v)\dd  v +\mathscr{R}_1,
\end{multline*}
where $\mathscr{R}_1$ is defined by
\begin{multline}\label{definition of remainder R_1}
\mathscr{R}_1(t,x):=2a^2m_g\int_{\R\times\RS}\M[\alpha\num,u,\theta/m_g](x,w)[F(x-a\n,v)-F(x,v)+a\n\cdot\nabla_xF(x,v)][(v-w)\cdot \n]^2\n\\
H((v-w)\cdot \n)\dd v\dd w\dd \n.
\end{multline}
Note that $\mathscr{R}_1$ is of order $o(a^3)$. 
We take the gradient of the volume fraction function $\alpha$ defined in \eqref{definition of volume fraction},
\[
\nabla_x\alpha=-\frac{4\pi}{3}a^3\int_{\R}\nabla_xF(v)\dd v.
\]
Hence, we have
\begin{multline}\label{computation:exchange of momentum}
    \lim_{\eta,\delta\to 0}\int_{\R}\E_1^\eta[f^{\eta,\delta},F^{\eta,\delta}](w)m_gw\dd w  =\pi a^2\alpha\num\theta\int_{\R} F (v) \q  \dd  v \\
    +\alpha \num\theta\nabla_x\alpha-\alpha\rho\int_{\R}Q (v-u)\nabla_xF(v)\dd  v 
    +\mathscr{R}_1,
\end{multline}
where $Q $ is defined in \eqref{definition of Q}. By the product rule, we have
\begin{equation}\label{product rule in momentum equation}
\alpha\rho\int_{\R}Q(v-u) \nabla_xF(v)\dd v=\Div_x\left(\alpha\rho\int_{\R}Q(v-u) F(v)\dd v\right)-\int_{\R}\Div_x[\alpha\rho Q(v-u)]F(v)\dd v.
\end{equation}
Note that the volume fraction function remains close to $1$, since the gas remains the dominant component even in the thick sprays regime. Consequently, the product $\alpha\nabla_x\alpha$ behaves similarly to $\nabla_x\alpha$ up to higher-order terms. More precisely, we expand the volume function $\alpha$ by definition \eqref{definition of volume fraction}:
\begin{multline}\label{approximation of density production}
\alpha\nabla_x\alpha=\left(1-\frac{4\pi}{3}a^3\int_{\R}F(v)\dd v\right)\left(-\frac{4\pi}{3}a^3\int_{\R}\nabla_xF(v)\dd v\right)=-\frac{4\pi}{3}a^3\int_{\R}\nabla_xF(v)\dd v+o(a^3)\\
=\nabla_x\alpha+o(a^3).
\end{multline}
And the difference is defined to be $\mathscr{R}_2$:
\begin{equation}\label{definition of remainder R_2}
    \mathscr{R}_2:=(1-\alpha)\num\theta\nabla_x\alpha=o(a^3)
\end{equation}
Then, we substitute \eqref{product rule in momentum equation} and \eqref{definition of remainder R_2} into \eqref{computation:exchange of momentum}. It shows 
\begin{multline}\label{exchange part of momentum equation of the gas phase}
   \lim_{\eta,\delta\to 0}\int_{\R}\E_1^\eta[f^{\eta,\delta},F^{\eta,\delta}](w)m_gw\dd w  =\pi a^2\alpha\num\theta\int_{\R} F (v) \q  \dd  v +\num\theta\nabla_x\alpha\\
   -\Div_x\left(\alpha\rho\int_{\R}Q(v-u) F(v)\dd v\right)+\int_{\R}\Div_x[\alpha\rho Q(v-u)]F(v)\dd v +\mathscr{R},
\end{multline}
where $\mathscr{R}:=\mathscr{R}_1+\mathscr{R}_2$.
Combining \eqref{transport part of momentum equation of the gas phase} and \eqref{exchange part of momentum equation of the gas phase}, we obtain the momentum equation for the gas phase
\begin{equation}\label{conservation of momentum for gas phase}
    \p_t (\alpha\rho u)+\nabla_x( \alpha\rho u^{\otimes 2})+ \alpha\nabla_x(\num\theta )=m_g\int_{\R}  D(v-u)F(v)\dd  v-\Div_x\left(\alpha\rho\int_{\R}Q(v-u) F(v)\dd v\right) +\mathscr{R}(t,x),
\end{equation}
where the microscopic friction force $D\equiv D[\alpha\rho,u,\theta](t,x,v)$ is defined by 
\[
D[\alpha\rho,u,\theta](v):=\pi a^2\alpha\num\frac{\theta}{m_g}\q+\Div_x (\alpha\num Q(v-u)).
\]

\subsection{Gas phase energy equation}
Again, we multiply gas the phase equation in \eqref{dimensionless kinetic equations} by $\frac{m_g|w|^2}{2}$, and then integrate it with respect to $w$, note the following identity
\begin{multline*}
\lim_{\eta,\delta\to 0}\int_{\R}\frac{|w|^2}{2}\p_t f^{\eta,\delta}\dd w+\int_{\R}\frac{|w|^2}{2}\nabla_x f^{\eta,\delta}\dd w\\
=\p_t \left(\int_{\R}\M[\alpha\num,u,\theta/m_g](w)\frac{|w|^2}{2}\dd w\right)+\Div_x\left(\int_{\R}\M[\alpha\num,u,\theta/m_g](w)\frac{|w|^2}{2}w\dd w\right).
\end{multline*}
We process the same steps: we split $w$ into $w-u$ and $u$, and the terms including odd orders of $w-u$ vanishes. Therefore,
\begin{multline*}
\p_t \left(\int_{\R}\M[\alpha\num,u,\theta/m_g](w)\frac{|w|^2}{2}\dd w\right)+\Div_x\left(\int_{\R}\M[\alpha\num,u,\theta/m_g](w)\frac{|w|^2}{2}w\dd w\right)\\
=\p_t \left(\int_{\R}\M[\alpha\num,u,\theta/m_g](w)\frac{|w-u|^2}{2}\dd w+\frac{|u|^2}{2}\int_{\R}\M[\alpha\num,u,\theta/m_g](w)\dd w\right)\\
+\Div_x\bigg[\left(\int_{\R}\M[\alpha\num,u,\theta/m_g](w)\frac{|w-u|^2}{2}\dd w\right)u+\left(\int_{\R}\M[\alpha\num,u,\theta/m_g](w)(w-u)^{\otimes 2}\dd w\right)u\\
+\left(\int_{\R}\M[\alpha\num,u,\theta/m_g](w)\dd w\right)\frac{|u|^2}{2}u\bigg].
\end{multline*}
It follows,
\begin{equation}\label{transport part of energy equation for gas phase}
\lim_{\eta,\delta\to 0}\int_{\R}\frac{m_g|w|^2}{2}\p_t f^{\eta,\delta}\dd w+\int_{\R}\frac{m_g|w|^2}{2}\nabla_x f^{\eta,\delta}\dd w=\p_t \left(\alpha\rho \left(\frac{|u|^2}{2}+\frac{3\theta}{2m_g} \right)\right)+\Div_x\left( \alpha\rho u\left(\frac{|u|^2}{2}+\frac{5\theta}{2m_g} \right)\right).
\end{equation}
We repeat the same approach as shown in the momentum equation, 
\begin{multline*}
    \int_{\R}\E_1^\eta[f^{\eta,\delta},F^{\eta,\delta}]\frac{|w|^2}{2}\dd w  \\
    =\frac{2a^2}{(1+\eta)^2}\int_{\R\times\RS} f^{\eta,\delta} (w) F^{\eta,\delta} (v) \big(( v -w) \cdot \n\big)^2(v+\eta w)\cdot \n H(( v -w) \cdot \n)  \dd  v \dd  \n \dd  w\\
  -\frac{2a^3}{(1+\eta)^2}\int_{\R\times\RS}  f^{\eta,\delta} (w)  \big(( v -w)\cdot \n\big)^2(v+\eta w)\cdot \n\nabla_x F^{\eta,\delta} (v)\cdot \n H(( v -w) \cdot \n)  \dd  v \dd  \n \dd  w+\mathscr{P}_1^{\eta,\delta},
\end{multline*}
where $\mathscr{P}_1^{\eta,\delta}$ is given by
\begin{multline*}
\mathscr{P}_1^{\eta,\delta}:=2a^2\int_{\R\times\RS} f^{\eta,\delta} (x,w) [F^{\eta,\delta} (x,v)-F^{\eta,\delta} (x-a\n,v)+a\n\cdot\nabla_xF^{\eta,\delta} (x,v)] \\
\big(( v -w) \cdot \n\big)^2(v+\eta w)\cdot \n H(( v -w) \cdot \n)  \dd  v \dd  \n \dd  w    
\end{multline*}

As $\eta,\delta\to 0$, we integrate over $\n\in \mathbb{S}^2$. For the same reason, 
\begin{multline*}
    \lim_{{\eta,\delta}\to 0}\int_{\R}\E_1^\eta[f^{\eta,\delta},F^{\eta,\delta}]\frac{m_g|w|^2}{2}\dd w  =\pi a^2 \alpha\num\theta\int_{\R}F (v) \q\cdot v \dd  v \\
    -\frac{4\pi}{3}a^3\alpha\num\theta\int_{\R}v\cdot\nabla_xF(v)\dd v
    -\alpha\rho\int_{\R}  Q(v-u) :v\otimes \nabla_xF(v)  \dd  v +\mathscr{P}_1,
\end{multline*}
where $\mathscr{P}_1$ is defined as
\begin{multline}\label{definition of remainder term P_1}
    \mathscr{P}_1:=2a^2m_g\int_{\R\times\RS} \M[\alpha\num,u,\theta/m_g](w) [F (x,v)-F (x-a\n,v)+a\n\cdot\nabla_xF (x,v)] 
\big(( v -w) \cdot \n\big)^2
\\v\cdot \n H(( v -w) \cdot \n)  \dd  v \dd  \n \dd  w  
\end{multline}
Thanks to the mass conservation of the particles, and the definition of the volume fraction,
\[
-\p_t\alpha+\frac{4\pi}{3}a^3\int_{\R}v\cdot \nabla_xF\dd v=\frac{4\pi}{3}a^3\left(\int_{\R}\p_tF\dd v+\int_{\R}v\cdot\nabla_xF\dd v\right)=\frac{4\pi}{3}a^3\int_{\R}\E_2^\eta[F,f]\dd v=0,
\]
the energy exchange becomes
\begin{multline}\label{computation: energy exchange}
    \lim_{{\eta,\delta}\to 0}\int_{\R}\E_1^\eta[f^{\eta,\delta},F^{\eta,\delta}]\frac{m_g|w|^2}{2}\dd w =\pi a^2 \alpha \num \theta\int_{\R}F (v) \q\cdot v \dd  v \\ -\alpha \rho\int_{\R}  Q(v-u) :v\otimes \nabla_xF(v)  \dd  v 
    -\alpha\num\theta\p_t\alpha+\mathscr{P}_1.
\end{multline}
Since $\alpha$ is close to 1, we can approximate $\alpha\p_t\alpha=(1+O(a^3))\p_t\alpha=\p_t\alpha+O(a^6)$. We define the remainder term $\mathscr{P}_2$ as
\begin{equation}\label{definition of remainder P_2}
    \mathscr{P}_2:=-(1-\alpha)\num\theta\p_t\alpha.
\end{equation}
Substituting \eqref{definition of remainder P_2} into \eqref{computation: energy exchange}, we obtain
\begin{multline*}
    \lim_{\eta,\delta\to 0}\int_{\R}\E_1^\eta[f^{\eta,\delta},F^{\eta,\delta}](w)\frac{m_g|w|^2}{2}\dd w  =\pi a^2\alpha\num\theta\int_{\R}F (v) \q\cdot v \dd  v  \\-\alpha\rho\int_{\R}  Q(v-u) :v\otimes \nabla_xF(v)  \dd  v -\num\theta\p_t\alpha+\mathscr{P},
\end{multline*}
where $\mathscr{P}=\mathscr{P}_1+\mathscr{P}_2$
Next, we move the gradient outside the integral,
\begin{equation*}
    \alpha\rho\int_{\R}  Q(v-u) :v\otimes \nabla_xF(v)  \dd  v=\Div_x\left(\alpha\rho\int_{\R}F(v)  Q(v-u)~ v \dd  v\right)-\int_{\R}F(v)  \Div_x(\alpha\rho Q(v-u))\cdot v \dd  v. 
\end{equation*}
Recall the definition of friction force $D$, we obtain
\begin{multline}\label{exchange part of energy equation of the gas phase}
    \lim_{\eta,\delta\to 0}\int_{\R}\E_1^\eta[f^{\eta,\delta},F^{\eta,\delta}](w)        \frac{m_g|w|^2}{2}\dd w  =m_g\int_{\R}F (v) D(v-u)\cdot v \dd  v  -\Div_x\left(\alpha\rho\int_{\R}F(v)  Q(v-u)~ v \dd  v\right)\\ -\num\theta\p_t\alpha+\mathscr{P}.
\end{multline}
Combining \eqref{transport part of energy equation for gas phase} and \eqref{exchange part of energy equation of the gas phase}, we obtain the energy equation for the gas phase
\begin{multline}\label{conservation of energy for gas phase}
   \p_t \left(\alpha\rho \left(\frac{|u|^2}{2}+\frac{3\theta}{2m_g} \right)\right)+\Div_x\left( \alpha\rho u\left(\frac{|u|^2}{2}+\frac{3\theta}{2m_g} \right)\right)+\Div_x(\alpha \num \theta u)+\num\theta\p_t\alpha\\
   =m_g\int_{\R}F(v)  D(v-u)\cdot v \dd  v  -\Div_x\left(\alpha\rho\int_{\R}F(v)  Q(v-u)~ v \dd  v\right) +\mathscr{P}.
\end{multline}


\section{Conclusion}\label{sec:conclusion}
To obtain the commonly used model, we combine equations \eqref{conservation of mass for gas phase}, \eqref{conservation of momentum for gas phase}, and \eqref{conservation of energy for gas phase}, and make the following conventions:
\begin{itemize}
    \item The gas phase is ideal gas: $p(t,x)=\num(t,x)\theta(t,x)$,
    \item Total energy per unit mass is: $E(t,x)=\frac{|u(t,x)|^2}{2}+\frac{3\theta(t,x)}{2m_g}$.
\end{itemize}
Readers may wonder whether it's reasonable to assume that the gas governed by Enskog equation is still behaves as an ideal gas. This is justified by the fact that we use the Enskog collision integral specifically to describe the collisions \emph{between particle and gas}, while the collision among gas molecules themselves are still governed by Boltzmann collision integral $\B[f,f]$. 

A natural question arises: can we generalize the cross-section $b_{pg}(z, \omega)$ within the Enskog–Boltzmann collision integral to account for more complex particle interactions beyond hard-sphere collisions? If we adopt such a generalized cross-section, the familiar coefficient $\frac{4\pi}{3}a^3$ (corresponds to the volume of sphere with radius $a$), which arises in the standard Enskog theory for hard spheres of radius $a$, is no longer exact. 
In the generalized setting, the coefficient in front of the pressure gradient term $\nabla_x p$ instead depends on the angular integrals

$$
\int_{\mathbb{S}^2} b_{pg}(z, \omega) \, \mathrm{d}\omega \quad \text{and} \quad \int_{\mathbb{S}^2} b_{pg}(z, \omega)(z \cdot \omega)^2 \, \mathrm{d}\omega,
$$
which capture how the effective collision geometry and scattering vary with orientation $\omega$ and relative position $z$. These integrals account for more general interaction mechanisms, including anisotropic or soft potentials, rather than assuming a fixed collision diameter.
Consequently, the parameter $a$, previously interpreted as the hard-sphere radius, should be replaced with an effective interaction range. This interaction range reflects the spatial extent over which particles influence each other and aligns conceptually with the van der Waals correction to real gas behavior, where the microscopic structure of intermolecular forces alters the macroscopic thermodynamic quantities. Thus, in models employing a generalized cross-section, the geometric interpretation of $a$ must be reconsidered in terms of the physical nature of the particle interactions. See, for instance, §74 in \cite{LANDAU5}.

Finally, we must emphasize that this result is still far from a rigorous derivation. Although the well-posedness for suitable initial datum is obtained in \cite{EH2023} recently, as of the time of this writing, very little is known about the existence, uniqueness, or qualitative behavior of solutions to the coupled Enskog-Boltzmann equations. The well-posedness of coupled kinetic equation remains largely open at the time of writing. Moreover, the reminder terms $\mathscr{P},\mathscr{Q}$ and $\mathscr{R}$ are highly intricate, involving not only the expansions on distribution functions, but also the density product between two species.

\section*{Acknowledgments}
The author would like to thank Laurent Desvillettes, and Fran\c{c}ois Golse for their insightful discussions and contributions to the development of this work. Special thanks are due to Fr\'ed\'erique Charles for her careful proofreading and constructive feedback.

\bibliographystyle{alpha}
\bibliography{thick_sprays}

\appendix
\newpage
\section{Proof of Proposition \ref{prop:weak formula for enskog collision integral}}\label{apd:proof of prop 3}
\begin{proof}
We do the change of variables $(v,w,\n)\mapsto(v',w',-\n)$ for the gain terms of $\int_{\R}\E_1^\eta[f,F]\varphi(x,w)\dd w$ and $\int_{\R}\E_2^\eta[F,f]\varphi(x,v)\dd v$, the properties \eqref{change the sign of n}-\eqref{specular reflection} show
\begin{multline*}
\int_{\R}\E_1^\eta[f,F]\varphi(x,w)\dd w=\int_{\R\times\RS} a^2  \big(\varphi(x,w')-\varphi(x,w)\big)f (x , w) F (x-a \n, v)
   ( v - w) \cdot \n H(( v - w) \cdot \n)\\
   \dd  \n \dd  v\dd w,
\end{multline*}
and
\begin{multline*}
\int_{\R}\E_2^\eta[F,f]\varphi(x,v)\dd v=\int_{\R\times\RS} a^2 \big(\varphi(x,v')-\varphi(x,v)\big) f (x + a \n, w) F (x, v)
   ( v - w) \cdot \n H(( v - w) \cdot \n) \\
   \dd  \n \dd  v\dd w.
\end{multline*}
Therefore, we have
\begin{multline*}
\eta\int_{\R} \E_1^\eta[f,F] \varphi(x,w)\dd w+\int_{\R}\E_2^\eta[F,f]\varphi(x,v)\dd v\\
=\int_{\RS\times\R}a^2 f(x+a\n,w)F(x,v)(v-w)\cdot \n H((v-w)\cdot \n)\bigg[\varphi(x,w')+\varphi(x+a\n,v')\\
-\eta\varphi(x,w)-\eta\varphi(x+a\n,v)\bigg]\dd \n \dd v\dd w  \\
 +\eta\int_{\RS\times\R}a^2  (v-w)\cdot \n H((v-w)\cdot \n)\bigg[(\varphi(x,w')-\varphi(x,w)\big)f (x , w) F (x-a \n, v)\\
 -(\varphi(x+a\n,w')-\varphi(x+a\n,w)\big)f (x+a\n , w) F (x, v)\bigg]\dd \n \dd v\dd w .
\end{multline*}
The fundamental theorem of calculus shows
\[
\Phi(x+a\n)-\Phi(x)=\int^a_0\frac{\dd}{\dd s}\Phi(x+s\n)\dd s=\int^a_0\n\cdot \nabla_x\Phi(x+s\n)\dd s.
\]
Since $\varphi$ is smooth, and furthermore if $f(\cdot,w)F(\cdot,w)$ is smooth, then the last line of the right-hand side can be rewritten as 
\begin{multline*}
    -\eta\Div_x\int_{\RS\times\R}\int^a_0a^2 \n(v-w)\cdot \n H((v-w)\cdot \n)(\varphi(x+s \n,w')-\varphi(x+s\n,w)\big)\\f (x+s \n , w) F (x-(a-s) \n, v)
    \dd \n \dd v\dd w\dd s .
\end{multline*}
\end{proof}

\newpage
\section{Two integrals}\label{apd:two integrals}
Let \( \xi \in \mathbb{R}^3\setminus\{0\} \) be a fixed vector. We compute these two integrals

\[
K(\xi):=\int_{\mathbb{S}^2} (\omega\cdot \xi)^2\omega H(\omega \cdot \xi) \, \mathrm{d}\omega
\quad \text{and} \quad
K_d(\xi):=\int_{\mathbb{S}^{d-1}} \omega^{\otimes 2}(\omega\cdot \xi)^2 H(\omega \cdot \xi) \, \mathrm{d}\omega.
\]
In the following, we adopt the Einstein summation convention for tensor calculation. $\delta_{ij}$ is the Kronecker delta. The colon symbol $:$ denotes a double contraction. For $2$-tensors $A\equiv (A_{ij})$ and $B\equiv (B_{ij})$,  We define $A:B=A_{ij}B_{ij}$.

\textbf{Compute $K(\xi)$}:

Because of the rotational symmetry of the sphere, we can assume without loss of generality that \( \frac{\xi}{|\xi|} = e_3 = (0,0,1) \), and express the unit vector in spherical coordinates as:

\[
\omega = (\sin\theta \cos\phi, \sin\theta \sin\phi, \cos\theta)
\]
with \( \theta \in [0, \pi] \), \( \phi \in [0, 2\pi) \), and the surface element:

\[
\mathrm{d}\omega = \sin\theta \, \mathrm{d}\theta \, \mathrm{d}\phi.
\]

The function $H(\omega\cdot\xi)$ becomes \( H(\cos\theta) \), which restricts the integral to the northern hemisphere \( \theta \in [0, \pi/2] \).
\[
K(\xi)=|\xi|^2\int_{0}^{2\pi}\int_0^{\pi/2}\omega \cos^2\theta\sin \theta\mathrm{d}\theta \, \mathrm{d}\phi.
\]
Noting that $\phi$ goes through $[0,2\pi]$, the entry containing $\sin\phi$ or $\cos\phi$ vanishes. Thus
\[
K(\xi)={2\pi}|\xi|^2\int_0^{\pi/2}e_3\cos^3\theta\sin \theta\mathrm{d}\theta \, \mathrm{d}\phi=\frac{\pi}{2}|\xi|\xi.
\]

\textbf{Compute $K_d(\xi)$}:
We identify $K_d(\xi)$ as the contraction of a $4$-tensor with a $2$-tensor:
\[
\left(\int_{\mathbb{S}^{d-1}}\omega^{\otimes 4}H(\omega\cdot\xi)\dd \omega\right):\xi^{\otimes 2}.
\]
By substituting $\omega\to -\omega$ in $\int_{\mathbb{S}^2}\omega^{\otimes 4}H(\omega\cdot\xi)\dd \omega$, we observe that
\[
\int_{\mathbb{S}^{d-1}}\omega^{\otimes 4}H(\omega\cdot\xi)\dd \omega=\frac{1}{2}\int_{\mathbb{S}^{d-1}}\omega^{\otimes 4}\dd \omega.
\]
This is a rank-$4$ isotropic and symmetric tensor, whose general form is given by:
\[
\int_{\mathbb{S}^{d-1}} \omega_i\omega_j\omega_k\omega_l \mathrm{d}\omega = C_d(\delta_{ij}\delta_{kl}+\delta_{ik}\delta_{jl}+\delta_{il}\delta_{jk}).
\]
To determine the constant $C_d$, we contract both sides by $\delta_{ij}\delta_{kl}$ :
\[
\delta_{ij}\delta_{kl}\int_{\mathbb{S}^{d-1}} \omega_i\omega_j\omega_k\omega_l \mathrm{d}\omega = C_d\delta_{ij}\delta_{kl}(\delta_{ij}\delta_{kl}+\delta_{ik}\delta_{jl}+\delta_{il}\delta_{jk}).
\]
It follows
\[
\int_{\mathbb{S}^{d-1}} (\omega_i\omega_i)(\omega_k\omega_k) \mathrm{d}\omega = C_d[(\delta_{ii})(\delta_{kk})+2\delta_{ik}\delta_{ik}].
\]
In the three dimensions, we obtain $4\pi=\int_{\mathbb{S}^2} 1 \mathrm{d}\omega = 15C_3$, namely, $C_3=\frac{4\pi}{15}$. Substituting back into $K_d(\xi)$, its entries are
\[
(K_3(\xi))_{ij}=\frac{2\pi}{15}(\delta_{ij}\delta_{kl}+\delta_{ik}\delta_{jl}+\delta_{il}\delta_{jk})\xi_k\xi_l=\frac{2\pi}{15}(\delta_{ij}(\xi_{k}\xi_k)+\xi_i\xi_j+xi_j\xi_i).
\]
It can be also written as
\[
K_3(\xi)=\frac{2\pi}{15}(|\xi|^2\id+2\xi^{\otimes 2}).
\]

\newpage
\section{Generalization to Dimension $d$}\label{apd:general dimension}
After taking limit, the coefficient in front of the pressure modification appears solely in 
\[
a^{d} \int_{\mathbb{R}^d\times\mathbb{S}^{d-1}} \n \cdot [\nabla_x \M (x, w') F^{\eta,\delta} (x, v') +
   \nabla_x \M (x, w) F (x, v)] ( v - w) \cdot \n  H(( v -w) \cdot \n)\dd  \n \dd  w
\]
It's not difficult to verify that, in $d$-dimensional space, it can be reacted as
\begin{equation*}
    2a^dF(v)\Div_x\int_{\mathbb{R}^d}K_d(v-w)\M(w)\dd w.
\end{equation*} 
To determinate $K_d$ which we've defined in Appendix \ref{apd:two integrals}, we note the identity
\[
|\mathbb{S}^{d-1}|=\int_{\mathbb{S}^{d-1}}1\dd \omega=C_d(d^2+2d).
\]
Therefore, 
\begin{equation*}
    \mathcal{R}_2[F,\M]=a^dF(v)\frac{|\mathbb{S}^{d-1}|}{d(d+2)}\Div_x\int_{\mathbb{R}^d}[2(v-w)^{\otimes 2}+|v-w|^2\id]\M\dd w.
\end{equation*} 
Then we decompose $v-w=(v-u)+(u-w)$, while the terms involving $v-u$ don't contribute to pressure correction. Using the identity:
\[
\int_{\mathbb{R}^d}(w-U)^{\otimes 2}\M[R,U,T]\dd w=RT\id,
\]
we obtain, in dimension $d$,
\[
a^d\frac{|\mathbb{S}^{d-1}|}{d(d+2)}F(v)\Div_x(2\alpha\num\theta\id+d\alpha\num\theta\id)=a^d\frac{|\mathbb{S}^{d-1}|}{d}F(v)\nabla_x(\alpha\num\theta).
\]
Since $d|B_d|=|\mathbb{S}^{d-1}|$, where $B_d$ denotes the unit ball in $\mathbb{R}^d$, we conclude that the coefficient in front of the pressure modification is exactly the volume of a ball of radius $a$ in $d$-dimensional space.

\newpage
\section{More on the friction force}\label{apd:more on the friction force}
As shown in \cite{DGR2019}, the leading-order friction force $q\left(\frac{v-u}{\sqrt{\theta/m_g}}\right)$ can be expressed as
\begin{equation}\label{leading-order friction force}
q(\xi)=\bar{q}(|\xi|)\xi
\end{equation}
Furthermore, we represent 
\[
\Div_x(\alpha \rho Q(v-u))=\frac{4\pi}{15}a^3\Div_x[2\alpha\rho(v-u)^{\otimes 2}+\alpha\rho|v-u|^2\id]
\]as a viscous tensor acting on relative velocity $v-u$ by using indices to represent the tensors to express the tensor component
\begin{multline*}
 \left(\Div_x \left(\alpha\rho(v-u)^{\otimes 2}\right)\right)_{j}=\p_{x_i}(\alpha\rho(v_i-u_i)(v_j-u_j))\\
 =\p_{x_i}(\alpha\rho)(v_i-u_i)(v_j-u_j)-\alpha \rho\p_{x_i}u_i(v_j-u_j)-\alpha\rho(v_i-u_i)\p_{x_i}u_j
\end{multline*}
hence
\begin{equation}\label{cpt:q_2bar1}
\Div_x \left(\alpha\rho(v-u)^{\otimes 2}\right)=\left[\nabla_x(\alpha\rho)\cdot(v-u)\id-\alpha\rho(\Div_x u)\id-\alpha\rho\nabla_xu\right](v-u).
\end{equation}
Similarly,
\[
 \left(\nabla_x \left(\alpha\rho|v-u|^2\right)\right)_{j}=\p_{x_j}(\alpha\rho(v_i-u_i)(v_i-u_i))=\p_{x_j}(\alpha\rho)(v_i-u_i)(v_i-u_i)-2\alpha\rho\p_{x_j}u_i(v_i-u_i),
\]
so that

\begin{equation}\label{cpt:q_2bar2}
\nabla_x \left(\alpha\rho|v-u|^2\right)=\left[\nabla_x(\alpha\rho)\otimes (v-u)-2\alpha\rho(\nabla_x u)^T\right](v-u).
\end{equation}
By substituting \eqref{cpt:q_2bar1} and \eqref{cpt:q_2bar2} into $\Div_x(\alpha\rho Q(v-u))$, and combining with \eqref{leading-order friction force}, we can the friction force  in the form of viscous tensor applying on relative velocity
\begin{equation}
D[\alpha\num,u,\theta]=\bar{D}[\alpha\num,v-u,\theta,u](v-u),
\end{equation}
where 
\begin{multline}\label{definition of barD}
\bar{D}[\alpha\num,v-u,\theta,u]:=\pi a^2 \alpha\num \sqrt{\frac{m_g}{\theta}}\bar{q}\left(\frac{|v-u|}{\sqrt{\theta/m_g}}\right)\id\\
+\frac{2\pi}{15}a^3\left[2(\nabla_x(\alpha\rho)\cdot(v-u)) \id+\nabla_x(\alpha\rho)\otimes (v-u)\right]
-\frac{4\pi}{15}a^3\alpha\rho\left[\Div_x u\id+\nabla_xu+(\nabla_xu)^T\right].
\end{multline}

Note that both $\pi a^2\alpha\num \sqrt{\frac{m_g}{\theta}}\, \bar{q}\left(\frac{|v-u|}{\sqrt{\theta/m_g}}\right)\id$ and $\frac{4\pi a^3}{15}(\nabla_x(\alpha\rho)\cdot(v-u)) \id$ are scalar multiples of the identity matrix and thus do not alter the direction of the force. In contrast, the term $\left[\nabla_x(\alpha\rho)\otimes (v-u)\right](v-u) = |v-u|^2 \nabla_x(\alpha\rho)$ introduces a directional deflection aligned with $\nabla_x(\alpha\rho)$. The expression $\Div_x u \, \id + \nabla_x u + (\nabla_x u)^T$ forms a symmetric tensor depending solely on the bulk velocity field of the gas phase.

\end{document}